\documentclass[reqno,a4paper]{amsart}
\usepackage{amssymb}    
\usepackage{latexsym}   

\newcommand {\IN}{\mathbb{N}}                          
\newcommand {\IZ}{\mathbb{Z}}
\newcommand {\IR}{\mathbb{R}}                           
\newcommand {\IC}{\mathbb{C}}
\newcommand {\T} {\mathbb T}

\newcommand   {\B}{\mathcal B}
\newcommand   {\D}{\mathcal D}
\newcommand   {\G}{\mathcal G}
\renewcommand {\H}{\mathcal H}
\newcommand {\I}{\mathcal I}
\newcommand   {\J}{\mathcal J}
\newcommand   {\K}{\mathcal K}
\renewcommand {\L}{\mathcal L}
\renewcommand {\P}{\mathcal P}

\newcommand   {\g}{\mathfrak{g}}
\newcommand   {\gc}{\mathfrak{g}_{_{\IC}}}
\renewcommand {\t}{\mathfrak{t}}
\newcommand   {\tc}{\t_{\IC}}                  
\newcommand   {\kk}{\mathfrak{K}}

\newcommand {\End}{\operatorname{End}}
\newcommand {\<}{\langle}
\renewcommand {\>}{\rangle}

\newcommand {\half}[1]{\frac{#1}{2}}

\newcommand {\Diff}{\operatorname{Diff}(S^{1})}
\newcommand {\Diffp}{\operatorname{Diff}_{+}(S^{1})}

\newcommand {\Ad}{\operatorname{Ad}}
\newcommand {\id}{\operatorname{id}}
\newcommand {\eexp}[2][]{e^{#1\overline{\pi(#2)}}}
\newcommand {\bbar}[1]{\overline{\pi(#1)}}
\newcommand {\Exp}{\operatorname{Exp}}
\newcommand {\hsmooth}{\mathcal H^{\infty}}

\newcommand {\der}[1]{\frac{d}{d#1}}
\newcommand {\diff}[2]{\left.\frac{d}{d#1}\right|_{#1=#2}}

\newcommand {\lldots}{\ldots\hskip -.15mm}
\newcommand {\wt}{\widetilde}
\newcommand {\wh}{\widehat}

\newcommand {\pback}{\rho^{*}U(\H)}
\newcommand {\phase}{\alpha_{\xi}}
\newcommand {\pol}{^{\scriptscriptstyle{\operatorname{pol}}}}
\newcommand {\lpol}{L\pol}

\newcommand {\vect}{\operatorname{Vect}(S^{1})}

\newcommand {\vectpolc}{\operatorname{Vect}_{\IC}\pol(S^{1})}
\newcommand {\Vir}{\operatorname{Vir}}

\newtheorem {theorem}{Theorem}[subsection]
\newtheorem {proposition}[theorem]{Proposition}
\newtheorem {corollary}[theorem]{Corollary}
\newtheorem {lemma}[theorem]{Lemma}

\newcommand {\ssection}[1]
            {\section{#1}\setcounter{theorem}{0}\setcounter{equation}{0}}
\newcommand {\ssubsection}[1]
            {\subsection{#1}\hfill\vspace{0.5\baselineskip}}
\renewcommand {\thetheorem}
{\ifcase\arabic{section}\else\arabic{section}.\ifcase\arabic{subsection}\else\arabic{subsection}.\fi\fi\arabic{theorem}}

\numberwithin {equation}{subsection}
\renewcommand {\theequation}
{\ifcase\arabic{section}\else\arabic{section}.\ifcase\arabic{subsection}\else\arabic{subsection}.\fi\fi\arabic{equation}}

\newcommand   {\ie}{{\it i.e. }}
\newcommand   {\eg}{{\it e.g. }}
\renewcommand {\proof}{{\sc Proof.}\ }
\newcommand   {\halmos}{$\Box$}
\newcommand   {\remark}{{\sc Remark.\ }}
\newcommand   {\definition}{{\sc Definition.\ }}

\begin{document}



\begin{center}
{\tt Journal of Functional Analysis {\bf 161} (1999), 478-508}
\end{center}

\title
[Integrating Unitary Representations of Lie Groups]
{Integrating Unitary Representations of Infinite--Dimensional Lie Groups}
\author{Valerio Toledano Laredo}
\address{
Institut de Math\'ematiques de Jussieu\\
UMR 7586\\
Case 191\\
Universit\'e Pierre et Marie Curie\\
4, Place Jussieu\\
F--75252 Paris Cedex 05}
\email{toledano@math.jussieu.fr}
\thanks{
Work supported by a TMR fellowship, contract no. FMBICT950083
and by a Non--Commutative Geometry Network grant, contract no.
ERB FMRXCT960073, both from the European Commission.}
\begin{abstract}
We show that in the presence of suitable commutator estimates,
a projective unitary representation of the Lie algebra of a
connected and simply connected Lie group $\G$ exponentiates to
$\G$. Our proof does not assume $\G$ to be finite--dimensional
or of Banach--Lie type and therefore encompasses the diffeomorphism
groups of compact manifolds.
We obtain as corollaries short proofs of Goodman and Wallach's
results on the integration of positive energy representations
of loop groups and Diff$(S^{1})$ and of Nelson's criterion for
the exponentiation of unitary representations of
finite--dimensional Lie algebras.
\end{abstract}
\maketitle

\ssection{Introduction}

The integration of unitary representations of a finite--dimensional Lie
group $\G$ from those of its Lie algebra has been well--understood
since the fundamental work of Nelson \cite{Ne1}. Using analytic vectors,
one formally regards the unitary group $U(\H)$
of the corresponding Hilbert space as an analytic Lie group and
obtains a local homomorphism $\G\longrightarrow U(\H)$ via the
Baker--Campbell--Hausdorff formula. This, and more recent methods
(see \eg \cite{Ro}), do not however apply when $\G$ is infinite--dimensional
for in the absence of a general inverse function theorem, its
exponential map may fail to be locally one--to--one, as is the
case for the diffeomorphism groups of compact manifolds \cite{Mi1}.\\

In the present paper, we describe a method to exponentiate a
unitary representation of the Lie algebra $\L$ of $\G$ when
the action of $\L$ is controlled by suitable commutator
estimates. Our method does not rely on the use of the
exponential map of $\G$ and applies to finite and
infinite--dimensional Lie groups, whether of Banach--Lie type
or not. It allows moreover to deal with projective representations.
More precisely, let $\pi:\L\longrightarrow\End(V)$ be a projective
representation of $\L$ by skew--symmetric operators acting on a dense
subspace $V$ of a Hilbert space $\H$. Thus, $\pi$ is linear and for
any $X,Y\in\L$
\begin{equation}
[\pi(X),\pi(Y)]=\pi([X,Y])+iB(X,Y)
\end{equation}
for some real--valued two--cocycle $B$ on $\L$.
Our main assumption is the existence of a self--adjoint operator
$A\geq 1$ on $\H$ for which $V$ is the space of smooth vectors,
\ie $V=\bigcap_{n\geq 0}\D(A^{n})$ and such that for any $\xi\in V$
and $n\in\IN$
\begin{align}
\|\pi(X)\xi\|_{n}&\leq |X|_{n+1}\|\xi\|_{n+1}
\label{assu 1}\\
\|[A,\pi(X)]\xi\|_{n}&\leq |X|_{A,n+1}\|\xi\|_{n+1}
\label{assu 2}
\end{align}
where $\|\xi\|_{n}=\|A^{n}\xi\|$ and the $|.|$ are continuous
semi--norms on $\L$.\\

A justification of this assumption might be in order. If $\G$ is
finite--dimensional, one usually assumes that the Laplacian
$\Delta=\sum_{i}\pi(X_{i})^{2}$ corresponding to some basis $X_{i}$
of $\L$ is essentially self--adjoint on $V$. Nelson's theorem then
guarantees that $\pi$ exponentiates to a unitary representation
of $\G$ \cite{Ne1}. In this case, the action of $\L$ extends to
$\bigcap_{n\geq 0}\D(\overline{\Delta}^{n})$ and, setting
$A=1-\overline{\Delta}$, the estimates \eqref{assu 1}--\eqref{assu 2}
follow from simple intrinsic manipulations in the enveloping algebra
of $\L$ \cite{Ne1} so that Nelson's criterion is a special case of
our assumptions. For an infinite--dimensional Lie group, no analogue
of the Laplacian
exists but the above estimates were noticed by Goodman and Wallach in
their work on positive energy representations of loop groups
$LG=C^{\infty}(S^{1},G)$ and of $\Diff$ \cite{GoWa1,GoWa2}. In both
cases, the 'laplacian' $A$ is the infinitesimal generator $L_{0}$ of
rotations which embed in $\Diff$ on the one hand and act automorphically
on $LG$ on the other. Unlike the finite--dimensional situation, the
estimates \eqref{assu 1}--\eqref{assu 2} are a consequence of
representation--dependent computations relying crucially on the
fact that $L_{0}$ has non--negative spectrum, a defining property
of positive energy representations.\\

The basic observation of the present paper is that the estimates
\eqref{assu 1}--\eqref{assu 2} may be used to regard the unitary group
$U(\H)$ as a {\it regular} Lie group. Formally then, a homomorphism
$\G\longrightarrow U(\H)$ can be constructed by using a variant of Lie
theory due to Thurston and Omori {\it et al.} \cite{Om,Mi1}. Recall
from \cite{Mi1} that a Lie group $\G'$ with Lie algebra $\L'$ is said
to be regular if for any $X\in C^{\infty}(I,\L')$, with $I=[0,1]$,
there exists $p\in C^{\infty}(I,\G')$ such that $\dot p=Xp$ and
$p(0)=1$. $p$ is then unique so that the product integral or Volterra
map $X\longrightarrow p(1)$ is well--defined and assumed to be smooth.
This is a time--dependent exponential map, for if
$X(t)\equiv X_{0}\in\L$, then $p(t)=\exp_{\G}(tX_{0})$.\\

Product integrals may be used as a substitute for 
Baker--Campbell--Hausdorff series as follows. Let $\G,\G'$ be Lie
groups with Lie algebras $\L,\L'$ and assume that $\G$ is connected
and simply connected and $\G'$ regular. Then, any continuous
homomorphism $F:\L\longrightarrow\L'$ determines a unique
homomorphism $\Phi:\G\longrightarrow\G'$ with differential $F$ as
follows. Let $g\in\G$ and $p$ a smooth path in $\G$ with $p(0)=0$
and $p(1)=g$. Let $X=\dot pp^{-1}\in C^{\infty}(I,\L)$ and
$q\in C^{\infty}(I,\G')$ be such that $\dot q=F(X)q$ and $q(0)=1$.
Set $\Phi(g)=q(1)$. To see that this is independent of $p$, pick
a smooth homotopy $H:I^{2}\rightarrow \G$ with $H(0,\cdot)
\equiv 1$, $H(1,\cdot)\equiv g$ and $H(t,0)=p(t)$. The partial
derivatives $\partial_{i}H H^{-1}$ define a flat $\G$--connection
on $I^{2}$. Composing with $F$, we get a flat $\G'$--connection 
a horizontal section $s$ of which may be constructed
using product integrals. Since the connection vanishes on $\{0,1\}
\times I$, we get $s(1,0)s(0,0)^{-1}=s(1,1)s(0,1)^{-1}$ so that
$\Phi$ is well--defined and is easily seen to be a
homomorphism.\\

In our representation--theoretic context where, formally
$\G'=U(\H)$, the rigorous definition of $\Phi:\G\longrightarrow U(\H)$
amounts to a study of the time--dependent Schr\"odinger equation
\begin{align}
\frac{d\xi(t)}{dt}&=\pi(X(t))\xi(t)
\label{eq:sch 1}\\
\xi(0)&=\xi_{0}
\label{eq:sch 2}
\end{align}
determined by $X\in C^{\infty}(I,\L)$ and $\xi_{0}\in V$ and most
of the paper is devoted to proving the smooth well--posedness of
\eqref{eq:sch 1}--\eqref{eq:sch 2}. Once that is established, the
required exponentiation of $\pi$ is obtained by sending $g\in\G$
to the unitary operator $U(g)$ mapping $\xi_{0}\in V$ to $\xi(1)$,
where $\xi$ is the unique solution of \eqref{eq:sch 1}--\eqref{eq:sch 2}
with $X=\dot pp^{-1}$ and $p$ a smooth path in $\G$ with $p(0)=1$
and $p(1)=g$.\\

The paper is structured as follows. In \S \ref{se:ODE}, we
consider the time--independent version of
\eqref{eq:sch 1}--\eqref{eq:sch 2}.
Using Nelson's commutator theorem, we prove that the action of $\L$
on $V$ is essentially skew--adjoint and that the corresponding
one--parameter groups preserve the scale defined by $A$.
In \S \ref{se:tODE}, we prove the continuous well--posedness of
\eqref{eq:sch 1}--\eqref{eq:sch 2} by using product integrals
with unbounded generators. The smooth well--posedness is established
in \S \ref{se:posedness} by studying the inhomogeneous equation
obtained by formally differentiating
\eqref{eq:sch 1}--\eqref{eq:sch 2} with respect to $X$.
\S \ref{se:exponentiating} contains our main result
(theorem \ref{th:main}). We define a Volterra map :
$C(I,\L)\longrightarrow U(\H)$ and show that it factors through
a projective unitary representation $\rho$ of $\G$ the differential
of which is $\pi$. We prove moreover that the central extension of
$\G$ induced by $\rho$ is smooth. Finally, in \S \ref{se:applications}
we apply our exponentiation result to positive energy representations
of loop groups and $\Diff$ and unitary representations of
finite--dimensional Lie algebras.\\

{\bf Acknowledgements.} I am grateful to Professor G. Segal for a
number of very stimulating conversations on infinite--dimensional
Lie groups and to Professor J. Eells for his interest in the present
work. This paper was begun at the Department of Pure Mathematics and
Mathematical Statistics of the University of Cambridge and
completed at the Institut de Math\'ematiques de Jussieu of
the Universit\'e Pierre et Marie Curie, within the Alg\`ebres
d'Op\'erateurs et Repr\'esentations group. I wish to thank both
institutions for their kind hospitality and pleasant working
atmospheres.

\ssection{Time--independent ODE's in $\H$}
\label{se:ODE}

Throughout this paper, $\G$ denotes a connected and simply connected
Lie group and $\L$ its Lie algebra. We follow Milnor's terminology and
consider a Lie group to be a smooth manifold modelled on a complete,
locally convex, real topological vector space of finite or infinite
dimension, and possessing a compatible group structure \cite{Mi1}.
Let $\pi:\L\longrightarrow\End(V)$ be a projective representation of
$\L$ by skew--symmetric operators acting on a dense subspace $V$ of
a Hilbert space $\H$. We henceforth assume the existence of a
self--adjoint operator $A\geq 1$ for which
$V=\bigcap_{n\geq 0}\D(A^{n})$ and such that, for any $X\in\L$,
$\xi\in V$ and $n\in\IN$
\begin{align}
\|\pi(X)\xi\|_{n}&\leq |X|_{n+1}\|\xi\|_{n+1}
\label{sobolev one}\\
\|[A,\pi(X)]\xi\|_{n}&\leq |X|_{A,n+1}\|\xi\|_{n+1}
\label{sobolev two}
\end{align}
where $\|\xi\|_{n}=\|A^{n}\xi\|$ and the $|.|$ are
continuous semi--norms on $\L$ which for convenience
we take to be increasing in $n$
\footnote{we have departed from the usual convention that
$\|\xi\|_{n}=\|A^{\half{n}}\xi\|$ which implies that $A$
is of order 2 with respect to its own scale. Since all
operators we shall be considering are of the same order
as $A$, we have preferred to normalise that order to 1.}.\\

For any $s\in\IR$, let $\H^{s}$ be the completion of $V$ with
respect to the inner product $(\xi,\eta)_{s}=(A^{s}\xi,A^{s}\eta)$
so that $A$ defines unitaries $\H^{s}\rightarrow\H^{s-1}$ and,
if $s\geq 0$, $\H^{s}=\D(A^{s})$. Let also $\H^{\infty}$ be
$V=\bigcap_{s}\H^{s}$ with the corresponding Fr\'echet topology.
Since $(\xi,\eta)=(A^{s}\xi,A^{-s}\eta)$, $\H^{s}$ is canonically
isomorphic to the (anti--)dual of $\H^{-s}$. In particular, by
skew--symmetry of the $\pi(X)$, the estimates
\eqref{sobolev one}--\eqref{sobolev two} extend to any $n\in\IZ$
provided we set $|X|_{-n}=|X|_{n+1}$ for $n\geq 0$.
By \eqref{sobolev one}, the operators $\pi(X)\in\End(\hsmooth)$
extend to bounded linear maps $\H^{s}\rightarrow\H^{s-1}$.
In particular, if $s\geq 1$, $\pi_{s}(X):=\left.\pi(X)\right|_{\H^{s}}$
are densely defined skew--symmetric and therefore closeable
operators on $\H$ which, by \eqref{sobolev one} have a common
closure $\bbar{X}$. We shall loosely refer to any of
the $\pi_{s}(X)$ as $\pi(X)$, drawing a distinction only between
these and $\bbar{X}$.

\begin{proposition}\label{invariance}
The operators $\pi(X)$, $X\in\L$, are essentially skew--adjoint
on $\H^{\infty}$ and any of the $\H^{n}$, $n\geq 1$. Moreover,
for each $n\in\IN$, the unitaries $e^{\bbar{X}}$
restrict to bounded linear maps $\H^{n}\rightarrow\H^{n}$ with
\begin{equation}\label{exponential estimate}
\|e^{\bbar{X}}\|_{\B(\H^{n})}\leq e^{2n|X|_{A,n}}
\end{equation}
and therefore define continuous automorphisms of $\hsmooth$.
\end{proposition}
\proof
The essential skew--adjointness claim follows from Nelson's commutator
theorem \cite[prop. 2]{Ne3} since, by
\eqref{sobolev one}--\eqref{sobolev two} and interpolation, $\pi(X)$
and $[A,\pi(X)]$ define bounded operators
$\H^{\half{1}}\rightarrow\H^{-\half{1}}$.
Let now $n\geq 1$ and $N=A^{2n}$. We will show that
$\H^{n}=\D(N^{\half{1}})$ is invariant under $e^{\bbar{X}}$
by using an elegant trick of Faris and Lavine \cite[thm. 2]{FL}.
We begin by establishing a simple quadratic form inequality. Let
$\eta\in\hsmooth$, then by
$[A^{2n},\pi(X)]=\sum_{k=0}^{2n-1}A^{k}[A,\pi(X)]A^{2n-1-k}$
and \eqref{sobolev two}, we get
\begin{equation}
\begin{split}
(\pi(X)\eta,N\eta)+(N\eta,\pi(X)\eta)
&=([N,\pi(X)]\eta,\eta)\\
&\leq\|N^{\half{1}}\eta\|\|N^{-\half{1}}[N,\pi(X)]\eta\|\\
&\leq\|N^{\half{1}}\eta\|\sum_{k=0}^{2n-1}
 \|A^{k-n}[A,\pi(X)]A^{2n-1-k}\eta\|\\
&\leq\|N^{\half{1}}\eta\|\sum_{k=0}^{2n-1}
 |X|_{A,k-n+1}\|A^{n}\eta\|\\
&\leq 2n|X|_{A,n}\|N^{\half{1}}\eta\|^{2}\\
\end{split}
\end{equation}
whence, by continuity, for any $\xi\in\D(N)$
\begin{equation}\label{quadratic}
(\pi(X)\eta,N\eta)+(N\eta,\pi(X)\eta)\leq
2n|X|_{A,n}(\eta,N\eta)
\end{equation}
Let now $\epsilon>0$ and $N_{\epsilon}=N(\epsilon N+1)^{-1}$,
a bounded self--adjoint operator and notice that, by the spectral
theorem,
$\D(N^{\half{1}})=
 \{\xi\in\H|\lim_{\epsilon\searrow 0}(\xi,N_{\epsilon}\xi)<\infty\}$.
Fix $\xi\in\D(\bbar{X})$ and let
$\xi_{t}=e^{t\bbar{X}}\xi$, then
\begin{equation}\label{derivative}
\frac{d}{dt}(\xi_{t},N_{\epsilon}\xi_{t})=
(\bbar{X}\xi_{t},N_{\epsilon}\xi_{t})+
(N_{\epsilon}\xi_{t},\bbar{X}\xi_{t})
\end{equation}
To rewrite this differently, consider
\begin{equation}\label{try}
(\bbar{X}(\epsilon N+1)^{-1}\xi_{t},N_{\epsilon}\xi_{t})+
(N_{\epsilon}\xi_{t},\bbar{X}(\epsilon N+1)^{-1}\xi_{t})
\end{equation}
Using $(\epsilon N+1)^{-1}=1-\epsilon N_{\epsilon}$ and the fact that
$(\epsilon N+1)^{-1}$ maps $\H$ into $\D(N)\subset\D(\bbar{X})$
so that
$\epsilon N_{\epsilon}\xi_{t}=
 \xi_{t}-(\epsilon N+1)^{-1}\xi_{t}\in\D(\bbar{X})$,
we may rewrite \eqref{try} as
\begin{equation}
(\bbar{X}\xi_{t},N_{\epsilon}\xi_{t})-
\epsilon(\bbar{X}N_{\epsilon}\xi_{t},N_{\epsilon}\xi_{t})+
(N_{\epsilon}\xi_{t},\bbar{X}\xi_{t})-
\epsilon(N_{\epsilon}\xi_{t},\bbar{X}N_{\epsilon}\xi_{t})
\end{equation}
which, by the skew--adjointness of $\bbar{X}$ is equal to
\eqref{derivative}. Therefore, using \eqref{quadratic}
\begin{equation}
\begin{split}
\frac{d}{dt}(\xi_{t},N_{\epsilon}\xi_{t})
&=
(\bbar{X}(\epsilon N+1)^{-1}\xi_{t},N(\epsilon N+1)^{-1}\xi_{t})+
(N(\epsilon N+1)^{-1}\xi_{t},\bbar{X}(\epsilon N+1)^{-1}\xi_{t})\\
&\leq 2n|X|_{A,n}
((\epsilon N+1)^{-1}\xi_{t},N(\epsilon N+1)^{-1}\xi_{t})\\
&\leq 2n|X|_{A,n}
(\xi_{t},N_{\epsilon}\xi_{t})
\end{split}
\end{equation}
Integrating this inequality, we find
$(e^{t\bbar{X}}\xi,N_{\epsilon}e^{t\bbar{X}}\xi)\leq 
  e^{2n |X|_{A,n}|t|(\xi,N_{\epsilon}\xi)}$
for any $\xi\in\D(\bbar{X})$ and therefore for any $\xi\in\H$.
Choosing now $\xi\in\D(N^{\half{1}})$ with $\|\xi\|_{n}=1$ and
letting $\epsilon\rightarrow 0$ we see that
$e^{\bbar{X}}\xi\in\D(N^{\half{1}})$
and $\|e^{\bbar{X}}\xi\|_{n}\leq e^{2n |X|_{A,n}}$
as claimed \halmos

\begin{corollary}
For any $X\in\L$, $\xi\in\hsmooth$ and $k\geq 1$, we have
\begin{equation}\label{Taylor}
e^{(t+h)\bbar{X}}\xi=
e^{t\bbar{X}}\xi+\cdots+
\frac{h^{k}}{k!}\pi(X)^{k}e^{t\bbar{X}}\xi+
R(h)
\end{equation}
where all terms are in $\H^{\infty}$ and $R(h)=o(h^{k})$ in each
$\|\cdot\|_{n}$ norm, \ie $\|R(h)\|_{n}h^{-k}\rightarrow 0$ as
$h\rightarrow 0$.
\end{corollary}
\proof
We have $\xi\in\cap_{n}\D(\pi(X)^{n})\subset C^{\infty}(\bbar{X})$
and consequently, by Taylor's theorem \eqref{Taylor} holds where
\begin{equation}
R(h)=\int_{0}^{h}du_{1}\cdots\int_{0}^{u_{k}}du_{k+1}
     \pi(X)^{k+1}e^{u_{k+1}\bbar{X}}\xi
\end{equation}
is in $\hsmooth$ since all other terms are. Moreover, by
\eqref{sobolev one} and \eqref{exponential estimate}
\begin{equation}
\|R(h)\|_{n}
\leq\frac{|h|^{k+1}}{(k+1)!}
(|X|_{n+k+1})^{k+1}e^{2(n+k+1)|hX|_{A,n+k+1}}\|\xi\|_{n+k+1}
=o(h^{k})
\end{equation}
\halmos

\begin{corollary}
For any $X,Y\in\L$ and $\xi\in\H^{n+1}$, we have
\begin{equation}\label{exp diff estimates}
\|e^{\bbar{X}}\xi-e^{\bbar{Y}}\xi\|_{n}\leq
|X-Y|_{n+1}e^{2(n+1)\max(|X|_{A,n+1},|Y|_{A,n+1})}\|\xi\|_{n+1}
\end{equation}
\end{corollary}
\proof
Let $\xi\in\hsmooth$ and set $F(t)=e^{-t\bbar{X}}e^{t\bbar{Y}}\xi$.
By \eqref{Taylor}, $F$ is differentiable and
$\dot F=e^{-t\bbar{X}}(\pi(Y)-\pi(X))e^{t\bbar{Y}}\xi$. Thus,
\begin{equation}
\begin{split}
\|e^{\bbar{X}}\xi-e^{\bbar{Y}}\xi\|_{n}
&=\|e^{\bbar{X}}\int_{0}^{1}\dot F(t)dt\|_{n}\\
&\leq\int_{0}^{1}\|e^{(1-t)
  \bbar{X}}(\pi(Y)-\pi(X))e^{t\bbar{Y}}\xi\|_{n}dt\\
&\leq e^{2n|X|_{A,n}}|Y-X|_{n+1}e^{2(n+1)|Y|_{A,n+1}}\|\xi\|_{n+1}
\end{split}
\end{equation}
\halmos


\ssection{Time--dependent ODE's in $\H$}
\label{se:tODE}

\ssubsection{Product integrals with unbounded generators}
\label{ss:volterra}

If $X\in C(\IR,\L)$ and $a<b\in\IR$, we define below the
product integral $\prod_{b\geq\tau\geq a}\Exp(X(\tau)d\tau)$
by adapting the presentation of \cite[\S I.2]{Ne2} where the
case of bounded infinitesimal generators is treated. Consider
first step functions $X:[a,b]\rightarrow\L$,
$X(\tau)=X_{j}$ for $\tau_{j}\geq\tau>\tau_{j-1}$ corresponding to
subdivisions $b=\tau_{n}>\tau_{n-1}>\cdots>\tau_{1}>\tau_{0}=a$
and set
\begin{equation}
\prod_{b\geq\tau\geq a}\Exp(X(\tau)d\tau)\xi=
\eexp[\Delta_{n}]{X_{n}}\cdots\eexp[\Delta_{1}]{X_{1}}\xi
\end{equation}
where $\Delta_{j}=\tau_{j}-\tau_{j-1}$. If $X$, $Y$ are two step
functions, which we may take as defined on a common subdivision,
the identity
$E_{n}\cdots E_{1}-F_{n}\cdots F_{1}=\sum_{k=1}^{n}
 E_{n}\cdots E_{k+1}(E_{k}-F_{k})F_{k-1}\cdots F_{1}$
and the estimates \eqref{exponential estimate} and
\eqref{exp diff estimates} imply that
\begin{equation}\label{difference estimate}
\begin{split}
 &\|\negthickspace
  \prod_{b\geq\tau\geq a}\negthickspace\Exp(X(\tau)d\tau)\xi-
  \negthickspace
  \prod_{b\geq\tau\geq a}\negthickspace\Exp(Y(\tau)d\tau)\xi\|_{r}\\
\leq&
(b-a)|X-Y|_{r+1}^{[a,b]}
e^{2(r+1)(b-a)\max(|X|_{A,r+1}^{[a,b]},|Y|_{A,r+1}^{[a,b]})}\|\xi\|_{r+1}
\end{split}
\end{equation}
where $|Z|_{k}^{[a,b]}=\sup_{t\in[a,b]}|Z(t)|_{k}$. We may therefore
define the product exponential as a bounded operator
$\H^{r+1}\rightarrow\H^{r}$ for any $X\in C([a,b],\L)$ by using a
sequence of approximating step functions. However, since
\begin{equation}
\|\eexp[\Delta_{n}]{X_{n}}\cdots\eexp[\Delta_{1}]{X_{1}}\xi\|_{r}\leq
e^{2r(b-a)|X|_{A,r}^{[a,b]}}\|\xi\|_{r}
\end{equation}
and $\H^{r+1}$ is dense in $\H^{r}$, the operator extends to a
bounded linear map $\H^{r}\rightarrow\H^{r}$ satisfying
\begin{equation}\label{norm estimate}
\|\prod_{b\geq\tau\geq a}\Exp(X(\tau)d\tau)\|_{\B(\H^{r})}\leq
e^{2r(b-a)|X|_{A,r}^{[a,b]}}
\end{equation}
and \eqref{difference estimate}.
Notice that product integrals are invertible operators, in fact
\begin{equation}\label{eq:inversion}
 \prod_{b\geq\tau\geq a}\Exp(X(\tau)d\tau)^{-1}=
 \prod_{b\geq\tau\geq a}\Exp(-\check X(\tau)d\tau)
\end{equation}
where $\check X(\tau)=X(a+b-\tau)$ so that we may set, for $a>b$,
\begin{equation}
 \prod_{b\geq\tau\geq a}\Exp(X(\tau)d\tau)=
 \prod_{a\geq\tau\geq b}\Exp(X(\tau)d\tau)^{-1}
\end{equation}
With this convention, the following semigroup property always holds
for any $a,b,c\in\IR$
\begin{equation}\label{eq:semigroup}
\prod_{c\geq\tau\geq b}\Exp(X(\tau)d\tau)
\prod_{b\geq\tau\geq a}\Exp(X(\tau)d\tau)=
\prod_{c\geq\tau\geq a}\Exp(X(\tau)d\tau)
\end{equation}

\begin{lemma}\label{le:sot cont}
If $X\in C(\IR,\L)$, the map
$t\longrightarrow\prod_{t\geq\tau\geq 0}\Exp(X(\tau)d\tau)\in\B(\H^{r})$
is strongly continuous for any $r\in\IN$.
\end{lemma}
\proof
Since the operators $\prod_{t\geq\tau\geq 0}\Exp(X(\tau)d\tau)$ are
locally uniformly bounded in $\B(\H^{r})$, it is sufficient to check
strong continuity on the dense set of vectors $\xi\in\H^{r+1}\subset
\H^{r}$. By \eqref{eq:semigroup}
\begin{equation}
\begin{split}
 &
\|\negthickspace\prod_{t+h\geq\tau\geq 0}
  \negthickspace\Exp(X(\tau)d\tau)\xi-
  \negthickspace\prod_{t\geq\tau\geq 0}
  \negthickspace\Exp(X(\tau)d\tau)\xi\|_{r}\\
=&
\|(\negthickspace\prod_{t+h\geq\tau\geq t}
   \negthickspace\Exp(X(\tau)d\tau)-1)
   \negthickspace\prod_{t\geq\tau\geq 0}
   \negthickspace\Exp(X(\tau)d\tau)\xi\|_{r}\\
\leq&
|h||X|_{r+1}^{[t,t+h]}e^{2(r+1)|hX|_{A,r+1}^{[t,t+h]}}
\|\negthickspace\prod_{t\geq\tau\geq 0}
  \negthickspace\Exp(X(\tau)d\tau)\xi\|_{r+1}
\end{split}
\end{equation}
which tends to zero as $h\rightarrow 0$ \halmos

\ssubsection{Time--dependent ODE's in $\H^{r}$}

Let $I\ni 0$ be a fixed compact interval, $r\in\IN$, and
consider the Banach space $C(I,\H^{r+1})\cap C^{1}(I,\H^{r})$ 
with norm $f\rightarrow\|f\|_{r+1}^{I}+\|\dot f\|_{r}^{I}$
where $\|g\|_{k}^{I}=\sup_{t\in I}\|g(t)\|_{k}$.

\begin{theorem}\label{th:homo eq}
Let $X\in C(I,\L)$, $\xi\in\H^{r+1}$ and set
\begin{equation}
\I(X,\xi)(t)=\prod_{t\geq\tau\geq 0}\Exp(X(\tau)d\tau)\xi
\end{equation}
Then, $\I(X,\xi)\in C(I,\H^{r+1})\cap C^{1}(I,\H^{r})$ and
is the unique solution of
\begin{align}
\dot \xi(t)&=\pi(X(t))\xi(t) \label{eq:homo eq}\\
     \xi(0)&=\xi             \label{eq:homo in}
\end{align}
Moreover, the map
$C(I,\L)\times\H^{r+1}\longrightarrow C(I,\H^{r+1})\cap C^{1}(I,\H^{r})$,
$(X,\xi)\longrightarrow \I(X,\xi)$ is continuous.
\end{theorem}
\proof
The uniqueness follows from the skew--symmetry of the $\pi(X)$ since
for any solution of \eqref{eq:homo eq},
\begin{equation}
\frac{d}{dt}(\xi(t),\xi(t))=
(\pi(X(t))\xi(t),\xi(t))+(\xi(t),\pi(X(t))\xi(t))=
0
\end{equation}
and therefore $\|\xi(t)\|\equiv\|\xi(0)\|$. Let now $\I=\I(X,\xi)$.
By lemma \ref{le:sot cont}, $\I\in C(I,\H^{r+1})$. Moreover, by
the semigroup property,
$\I(t+h)-\I(t)=\prod_{t+h\geq\tau\geq t}\Exp(X(\tau)d\tau)\eta-\eta$
where $\eta=\I(t)\in\H^{r+1}$. By \eqref{difference estimate},
\begin{equation}
\|\negthickspace\prod_{t+h\geq\tau\geq t}
  \negthickspace\Exp(X(\tau)d\tau)\eta-
  \negthickspace\prod_{t+h\geq\tau\geq t}
  \negthickspace\Exp(X(t)d\tau)\eta\|_{r}
\leq
 |h||X-X(t)|_{r+1}^{[t,t+h]}e^{2(r+1)|X|_{A,r+1}^{[t,t+h]}}
\|\eta\|_{r+1}
=o(h)
\end{equation}
by continuity of $X$. On the other hand, in $\H^{r}$
\begin{equation}
\prod_{t+h\geq\tau\geq t}\Exp(X(t)d\tau)\eta=
\eexp[h]{X(t)}\eta=
\eta+h\pi(X(t))\eta+o(h)
\end{equation}
since $\eta\in\H^{r+1}\subset\D(\bbar{X(t)}^{r+1})$ so that
$\I(X,\xi)$ lies in $C^{1}(I,\H^{r})$ and satisfies
\eqref{eq:homo eq}--\eqref{eq:homo in}.\\

To check the continuity of $\I(X,\xi)$ in the first variable for
the $C(I,\H^{r+1})$ norm, let $X,Y\in C(I,\L)$, $\xi\in\H^{r+1}$
and $\eta\in\H^{r+2}$ be an auxiliary vector. Then,
\begin{equation}
\begin{split}
\|\I(X,\xi)-\I(Y,\xi)\|_{r+1}^{I}
&\leq
\|\I(X,\xi)-\I(X,\eta)\|_{r+1}^{I}+
\|\I(X,\eta)-\I(Y,\eta)\|_{r+1}^{I}+
\|\I(Y,\eta)-\I(Y,\xi)\|_{r+1}^{I}\\
&\leq
(e^{2(r+1)|I||X|_{A,r+1}^{I}}+e^{2(r+1)|I||Y|_{A,r+1}^{I}})
\|\xi-\eta\|_{r+1}\\
&+
|I||X-Y|_{r+2}^{I}e^{2(r+2)|I|\max(|X|_{A,r+2}^{I},|Y|_{A,r+2}^{I})}
\|\eta\|_{r+2}\\
\end{split}
\end{equation}
so that
\begin{equation}
\limsup_{Y\rightarrow X}\|\I(X,\xi)-\I(Y,\xi)\|_{r+1}^{I}
\leq
\inf_{\eta\in\H^{r+2}}2e^{2(r+1)|I||X|_{A,r+1}^{I}}\|\xi-\eta\|_{r+1}
=0
\end{equation}
Joint continuity in the $C(I,\H^{r+1})$ norm now follows from
\begin{equation}
\begin{split}
\|\I(X,\xi)-\I(Y,\psi)\|_{r+1}^{I}
&\leq
\|\I(X,\xi)-\I(Y,\xi)\|_{r+1}^{I}+
\|\I(Y,\xi)-\I(Y,\psi)\|_{r+1}^{I}\\
&\leq
\|\I(X,\xi)-\I(Y,\xi)\|_{r+1}^{I}+
e^{2(r+1)|I||Y|_{A,r+1}^{I}}\|\xi-\psi\|_{r+1}
\end{split}
\end{equation}
and in the $C^{1}(I,\H^{r})$ from
\begin{equation}
\begin{split}
\|\dot \I(X,\xi)-\dot \I(Y,\eta)\|_{r}^{I}
&=
\|\pi(X)\I(X,\xi)-\pi(Y)\I(Y,\eta)\|_{r}^{I} \\
&\leq
\|\pi(X)\I(X,\xi)-\pi(Y)\I(X,\xi)\|_{r}^{I}+
\|\pi(Y)\I(X,\xi)-\pi(Y)\I(Y,\eta)\|_{r}^{I}\\
&\leq
|X-Y|_{r+1}^{I}\|\I(X,\xi)\|_{r+1}^{I}+
  |Y|_{r+1}^{I}\|\I(X,\xi)-\I(Y,\eta)\|_{r+1}^{I}\\
\end{split}
\end{equation}
\halmos

\begin{corollary}\label{co:soI cont}
The map $C([a,b],\L)\rightarrow\B(\H^{r})$, 
$X\rightarrow\prod_{b\geq\tau\geq a}\Exp(X(\tau)d\tau)$
is strongly continuous for any $r\geq 0$.
\end{corollary}
\proof
For $r\geq 1$ the claim is a direct consequence of theorem
\ref{th:homo eq}. If $r=0$, the operators are unitaries in
$\B(\H)$ and it is sufficient to check strong continuity on
the dense subspace $\H^{1}\subset\H$. This in turn follows
from theorem \ref{th:homo eq} \halmos

\begin{corollary}\label{change of variable}
If $\phi:[a,b]\rightarrow[c,d]$ is a smooth map with $\phi(a)=c$
and $\phi(b)=d$, then, for any $X\in C([c,d],\L)$,
\begin{equation}\label{eq:change of variable}
\prod_{c\geq\tau\geq d}\Exp(X(\tau)d\tau)=
\prod_{a\geq\sigma\geq b}
\Exp(\phi'\cdot X\circ\phi(\sigma)d\sigma)
\end{equation}
\end{corollary}
\proof
Let $\xi_{0}\in\H^{1}$ and $\xi:[c,d]\rightarrow\H^{1}$ the
solution of $\dot\xi=\pi(X)\xi$, $\xi(c)=\xi_{0}$ so that
$\eta=\xi\circ\phi$ satisfies
$\dot\eta=\pi(\phi'\cdot X\circ\phi)\eta$ and $\eta(a)=\xi_{0}$.
Then, by uniqueness,
\begin{equation}
\prod_{c\geq\tau\geq d}\Exp(X(\tau)d\tau)\xi_{0}=
\xi(d)=
\eta(b)=
\prod_{a\geq\sigma\geq b}
\Exp(\phi'\cdot X\circ\phi(\sigma)d\sigma)\xi_{0}
\end{equation}
\halmos

\ssection{Smooth well--posedness of time--dependent ODE's in $\H$}
\label{se:posedness}

\ssubsection{Inhomogeneous differential equations in $\H^{r}$}

We shall be concerned with the continuous well--posedness of
the inhomogeneous linear equation
\begin{equation}
\dot\zeta=\pi(X)\zeta+\eta
\end{equation}
in $\zeta\in C(I,\H^{r+1})\cap C^{1}(I,\H^{r})$ with initial
condition $\zeta(0)=0$, $X\in C(I,\L)$ and $\eta\in C(I,\H^{r+1})$
which is obtained by formally differentiating
\eqref{eq:homo eq}--\eqref{eq:homo in}
with respect to $X$. A solution is
readily obtained through a variation of constants by setting
$\zeta(t)=\prod_{t\geq\tau\geq 0}\Exp(X(\tau)d\tau)\zeta_{0}(t)$
which yields
$\dot\zeta_{0}(t)= 
 \prod_{t\geq\tau\geq 0}\Exp(X(\tau)d\tau)^{-1}\eta(t)$.

\begin{theorem}\label{th:inhomo eq}
Let $X\in C(I,\L)$ and $\eta\in C(I,\H^{r+1})$, $r\geq 0$ and
define
\begin{equation}\label{inho sol}
\J(X,\eta)(t)=
\prod_{t\geq\tau\geq 0}\Exp(X(\tau)d\tau)
\int_{0}^{t}\prod_{s\geq\tau\geq 0}\Exp(X(\tau)d\tau)^{-1}\eta(s)ds=
\int_{0}^{t}\prod_{t\geq\tau\geq s}\Exp(X(\tau)d\tau)\eta(s)ds
\end{equation}
Then
\begin{enumerate}
\item[(i)]
$\J(X,\eta)\in C(I,\H^{r+1})\cap C^{1}(I,\H^{r})$ and is the
unique solution of
\begin{align}
\dot\J(X,\eta)   &=\pi(X)\J(X,\eta)+\eta \label{eq:inhomo eq}\\
    \J(X,\eta)(0)&=0                     \label{eq:inhomo in}
\end{align}
\item[(ii)] The map
$\J:C(I,\L)\times C(I,\H^{r+1})\longrightarrow
    C(I,\H^{r+1})\cap C^{1}(I,\H^{r})$, 
$(X,\eta)\longrightarrow \J(X,\eta)$ is continuous.
\end{enumerate}
\end{theorem}
\proof
(i) Uniqueness follows from that of solutions of the corresponding
homogeneous equation. $\J(X,\eta)\in C(I,\H^{r+1})$ because
$t\rightarrow\prod_{t\geq\tau\geq 0}\Exp(X(\tau)d\tau)\in\B(\H^{r+1})$
is strongly continuous and so is
\begin{equation}
s\longrightarrow\prod_{s\geq\tau\geq 0}\Exp(X(\tau)d\tau)^{-1}=
                \prod_{s\geq\tau\geq 0}\Exp(-X(s-\tau)d\tau)
\end{equation}
To check the differentiability of $\J=\J(X,\eta)$ in $\H^{r}$, write
\begin{equation}
\begin{split}
\J(t+h)-\J(t)
&=
\Bigl(
\negthickspace\prod_{t+h\geq\tau\geq t}
\negthickspace\Exp(X(\tau)d\tau)-1\Bigr)\J(t)\\
&+
\negthickspace\prod_{t+h\geq\tau\geq 0}
\negthickspace\Exp(X(\tau)d\tau)
\int_{t}^{t+h}
\negthickspace\prod_{s\geq\tau\geq 0}
\negthickspace\Exp(X(\tau)d\tau)^{-1}\eta(s)ds\\[1em]
&=
h\pi(X(t))\J(t)+o_{r}(h)\\[1em]
&+
\negthickspace\prod_{t+h\geq\tau\geq 0}
\negthickspace\Exp(X(\tau)d\tau)
(h\negthickspace\prod_{t\geq\tau\geq 0}
\negthickspace\Exp(X(\tau)d\tau)^{-1}\eta(t)+
o_{r+1}(h))\\[1em]
&=h\pi(X(t))\J(t)+h\eta(t)+o_{r}(h)\\
\end{split}
\end{equation}
where the subscript in $o_{k}(h)$ refers to the norm $\|\cdot\|_{k}$.\\
(ii) Let $X,Y\in C(I,\L)$ and $\eta,\psi\in C(I,\H^{r+1})$, then
\begin{equation}
\begin{split}
\|\J(X,\eta)(t)-\J(Y,\psi)(t)\|_{r+1}
&\leq
\int_{0}^{t}
\|\prod_{t\geq\tau\geq s}\Exp(X(\tau)d\tau)\eta(s)-
  \prod_{t\geq\tau\geq s}\Exp(Y(\tau)d\tau)\eta(s)\|_{r+1}ds\\
&+
\int_{0}^{t}
\|\prod_{t\geq\tau\geq s}\Exp(Y(\tau)d\tau)\eta(s)-
  \prod_{t\geq\tau\geq s}\Exp(Y(\tau)d\tau)\psi(s)\|_{r+1}ds\\
\end{split}
\end{equation}
The second term is bounded by
$|I|e^{2(r+1)|I||Y|_{A,r+1}^{I}}\|\eta-\psi\|_{r+1}^{I}$ and
tends to zero uniformly in $t$ as $Y\rightarrow X$,
$\psi\rightarrow\eta$. If $\xi\in C(I,\H^{r+2})$ is an auxiliary
function the first term is bounded by
\begin{equation}
\begin{split}
&\int_{0}^{t}
\|\prod_{t\geq\tau\geq s}\Exp(X(\tau)d\tau)\eta(s)-
  \prod_{t\geq\tau\geq s}\Exp(X(\tau)d\tau)\xi(s)\|_{r+1}ds\\
+&\int_{0}^{t}
\|\prod_{t\geq\tau\geq s}\Exp(X(\tau)d\tau)\xi(s)-
  \prod_{t\geq\tau\geq s}\Exp(Y(\tau)d\tau)\xi(s)\|_{r+1}ds\\
+&\int_{0}^{t}
\|\prod_{t\geq\tau\geq s}\Exp(Y(\tau)d\tau)\xi(s)-
  \prod_{t\geq\tau\geq s}\Exp(Y(\tau)d\tau)\eta(s)\|_{r+1}ds\\
\leq&
|I|e^{2(r+1)|I||X|_{A,r+1}^{I}}\|\eta-\xi\|_{r+1}^{I}\\
+&
|I|^{2}|X-Y|_{r+2}^{I}
e^{2(r+2)|I|\max(|X|_{A,r+2}^{I},|Y|_{A,r+2}^{I})}
\|\xi\|_{r+2}^{I}\\
+&
|I|e^{2(r+1)|I||Y|_{A,r+1}^{I}}\|\eta-\xi\|_{r+1}^{I}\\
\end{split}
\end{equation}
whence
\begin{equation}
\limsup_{\substack{Y\rightarrow X\\\psi\rightarrow\eta}}
\|\J(X,\eta)-\J(Y,\psi)\|_{r+1}^{I}
\leq\inf_{\xi\in C(I,\H^{r+2})}
2|I|e^{2(r+)|I||X|_{A,r+1}^{I}}\|\eta-\xi\|_{r+1}^{I}=0
\end{equation}
by density of the inclusion $C(I,\H^{r+2})\subset C(I,\H^{r+1})$.
The continuity of $\J$ in the $C^{1}(I,\H^{r})$ norm follows
easily from the above and the fact that
$\dot \J(X,\eta)=\pi(X)\J(X,\eta)+\eta$ \halmos\\

\remark Similar results hold if \eqref{eq:inhomo in} is replaced
by the initial condition $\zeta(0)=\zeta_{0}\in\H^{r}$ since the
solution is then given by
\begin{equation}
\zeta(t)=
\prod_{t\geq\tau\geq 0}\Exp(X(\tau)d\tau)\zeta_{0}+
\int_{0}^{t}\prod_{t\geq\tau\geq s}\Exp(X(\tau)d\tau)\eta(s)ds
\end{equation}
We won't however need to work in such generality.

\ssubsection{Differentiability properties of product exponentials}

We investigate below the smoothness of the map
$C(I,\L)\times\H\rightarrow\H$,
$(X,\xi)\rightarrow\prod_{1\geq\tau\geq 0}\Exp(X(\tau)d\tau)\xi$.

\begin{proposition}\label{pr:gateaux}
For any $X\in C(I,\L)$ and $\xi\in\H^{r+1}$, $r\in\IN$,
let $\I(X,\xi)\in C(I,\H^{r+1})\cap C^{1}(I,\H^{r})$ be
the unique solution of
\begin{equation}
\dot\I(X,\xi)=\pi(X)\I(X,\xi)
\end{equation}
with $\I(X,\xi)(0)=\xi$. Then, for any $1\leq m\leq r$ and
$\delta_{m},\lldots,\delta_{1}\in C(I,\L)$, the G\^ateaux
derivatives 
\begin{equation}
\I^{(m)}(X,\xi;\delta_{m},\lldots,\delta_{1})=
\lim_{h\rightarrow 0}
\frac{\I^{(m-1)}(X+h\delta_{m},\xi;\delta_{m-1},\lldots,\delta_{1})-
      \I^{(m-1)}(X,\xi;\delta_{m-1},\lldots,\delta_{1})}{h}
\end{equation}
with $I^{(0)}(X,\xi)=\I(X,\xi)$, exist in
$C(I,\H^{r-m+1})\cap C^{1}(I,\H^{r-m})$ and are continuous
in $(X,\xi,\delta_{m},\lldots,\delta_{1})$.
\end{proposition}
\proof
We shall prove inductively that
$I^{(m)}(X,\xi;\delta_{m},\lldots,\delta_{1})$ exists
in $C(I,\H^{r-m+1})\cap C^{1}(I,\H^{r-m})$, satisfies
the inhomogeneous differential equation
\begin{align}
\dot \I^{(m)}(X,\xi;\delta_{m},\lldots,\delta_{1})&=
\pi(X)\I^{(m)}(X,\xi;\delta_{m},\lldots,\delta_{1})+
\sum_{i=1}^{m}\pi(\delta_{i})
\I^{(m-1)}(X,\xi;\delta_{m},\lldots,\widehat\delta_{i},\lldots,\delta_{1})
\label{em}\\
\I^{(m)}(X,\xi;\delta_{m},\lldots,\delta_{1})(0)&=
\delta_{m,0}\xi
\label{iem}
\end{align}
and depends jointly continuously on $X,\delta_{m},\lldots,\delta_{1},\xi$.
For $m=0$, this follows from theorem \ref{th:homo eq}. Let now
$m\geq 1$ and set, for any $h\neq 0$,
\begin{equation}
\psi_{h}=\frac{1}{h}\Bigl(
         \I^{(m-1)}(X+h\delta_{m},\xi;\delta_{m-1},\lldots,\delta_{1})-
         \I^{(m-1)}(X,\xi;\delta_{m-1},\lldots,\delta_{1})\Bigr)
\end{equation}
By induction $\psi_{h}$ satisfies $\dot\psi_{h}=\pi(X)\psi_{h}+\eta_{h}$,
$\psi_{h}(0)=0$ where
\begin{equation}
\begin{split}
\eta_{h}
&=\pi(\delta_{m})
  \I^{(m-1)}(X+h\delta_{m},\xi;\delta_{m-1},\lldots,\delta_{1})\\
&+\sum_{i=1}^{m-1}\pi(\delta_{i})\frac{
\I^{(m-2)}
(X+h\delta_{m},\xi;\delta_{m-1},\lldots,\widehat\delta_{i},\lldots,\delta_{1})-
\I^{(m-2)}
(X,\xi;\delta_{m-1},\lldots,\widehat\delta_{i},\lldots,\delta_{1})}{h}
\end{split}
\end{equation}
Since the operators $\pi(\delta_{m}),\lldots,\pi(\delta_{1})$ define
bounded maps
$C(I,\H^{r-m+2})\longrightarrow C(I,\H^{r-m+1})$, we see by induction
that, as $h\rightarrow 0$, $\eta_{h}$ tends in $C(I,\H^{r-m+1})$
to the inhomogeneous term of \eqref{em}. It follows from theorem
\ref{th:inhomo eq} that $\psi_{h}$ tends in
$C(I,\H^{r-m+1})\cap C^{1}(I,\H^{r-m})$ to the unique solution
of \eqref{em}--\eqref{iem} and that the limit is continuous in
$X,\delta_{m},\lldots,\delta_{1},\xi$ \halmos

\begin{corollary}\label{smoothness}
The map $\I:C(I,\L)\times\hsmooth\longrightarrow\hsmooth$,
$(X,\xi)\longrightarrow\prod_{1\geq\tau\geq 0}\Exp(X(\tau)d\tau)\xi$
is smooth.
\end{corollary}
\proof
By proposition \ref{pr:gateaux}, all partial derivatives of $\I$
with respect to the first variable exist and are continuous. Since
$\I$ and these are linear in the second variable we deduce that $\I$
has continuous mixed partial derivatives of all orders and hence is
smooth \halmos

\begin{corollary}\label{smooth vectors}
The map $C(I,\L)\times\H^{r+1}\longrightarrow\H$, $r=0\ldots\infty$,
$(X,\xi)\longrightarrow\prod_{1\geq\tau\geq 0}\Exp(X(\tau)d\tau)\xi$
is of class $C^{r}$. In particular, any $\xi\in\H^{r+1}$ is of class
$C^{r}$ for the unitary action of $C(I,\L)$ on $\H$.
\end{corollary}

\remark For any $h\in\IR$ and $\xi_{0}\in\H^{r+1}$, let
$\xi_{h}\in C^{0}(I,\H^{r+1})\cap C^{1}(I,\H^{r})$ be the solution
of $\dot\xi_{h}=h\pi(X)\xi$, $\xi_{h}(0)=\xi_{0}$. Then, in $\H$
\begin{equation}
\begin{split}
\xi_{h}(t)
&=
\xi_{0}+h\int_{0}^{t}\pi(X(t_{1}))\xi_{h}(t_{1})dt_{1}\\
&=
\xi_{0}+
\sum_{k=1}^{r+1}h^{k}
\idotsint\limits_{t\geq t_{1}\geq\cdots\geq t_{k}\geq 0}
\pi(X(t_{1}))\cdots\pi(X(t_{k}))\xi_{0}dt_{1}\cdots dt_{k}+R(h)
\end{split}
\end{equation}
where
\begin{equation}
R(h)=
h^{r+1}\idotsint\limits_{t\geq t_{1}\geq\cdots\geq t_{r+1}\geq 0}
\pi(X(t_{1}))\cdots\pi(X(t_{r+1}))(\xi_{h}(t_{r+1})-\xi_{0})
dt_{1}\cdots dt_{r+1}=o(h^{r+1})
\end{equation}
since, as $h\rightarrow 0$, $\xi_{h}-\xi_{0}$ tends to zero in
$C^{0}(I,\H^{r+1})$.
Thus, if $\xi_{0}\in\hsmooth$ and $X\in C(I,\L)$ are fixed, the
Taylor series of the $C^{\infty}(\IR,\H)$ function
$h\longrightarrow\prod_{1\geq\tau\geq 0}\Exp(hX(\tau)d\tau)\xi_{0}$
at $h=0$ is given by the Dyson expansion
\begin{equation}
\xi_{0}+
\sum_{k\geq 1}h^{k}\idotsint\limits_{1\geq t_{1}\geq\cdots\geq t_{k}\geq 0}
\pi(X(t_{1}))\cdots\pi(X(t_{k}))\xi_{0}dt_{1}\cdots dt_{k}
\end{equation}

\remark Notice that the results and proofs of sections
\ref{se:ODE}--\ref{se:posedness} only depend on the fact that
$\L$ is a vector space and $\pi(\L)\longrightarrow\End(V)$ is
linear, not on the fact that $\L$ is a Lie algebra or $\pi$ a
representation.

\ssection{Integrating group representations}
\label{se:exponentiating}

\ssubsection{Uniqueness of exponentiation}
\label{ss:uniqueness of exp}

Let $\pi:\L\longrightarrow\End(\D)$ be a projective representation
of $\L$ by skew--symmetric operators acting on a dense subspace of
a Hilbert space $\H$ and $B$ the corresponding cocycle so that, for
any $X,Y\in\L$
\begin{equation}
[\pi(X),\pi(Y)]=\pi([X,Y])+iB(X,Y)
\end{equation}

\definition
An exponentiation of $\pi$ is a strongly continuous homomorphism
\begin{equation}
\rho:\G\longrightarrow PU(\H)=U(\H)/\T
\end{equation}
admitting an invariant subspace $\D_{\rho}$ with
\begin{equation}
\D\subseteq\D_{\rho}\subseteq\bigcap_{X\in\L}\D(\bbar{X})
\end{equation}
and such that for any $p\in C^{\infty}(\IR,\G)$ with $p(0)=1$,
there exists, for small $t$, a continuous lift $\wt\rho(p(t))$
of $\rho(p(t))$ such that for any $\xi\in\D_{\rho}$,
$t\rightarrow\wt\rho(p(t))\xi$ is of class $C^{1}$ and satisfies
\begin{equation}\label{eq:local ode}
\frac{d}{dt}\wt\rho(p)\xi=
\bbar{\dot{p}p^{-1}}\wt\rho(p)\xi
\end{equation}
If $\pi$ is an ordinary representation, \ie $B=0$, we demand
in addition that $\rho$ map into $U(\H)$ and that $\wt\rho=\rho$.\\

\remark The above definition is somewhat stronger than the usual
ones but avoids the use of one--parameter groups which may fail
to exist or to generate $\G$ if the latter is an arbitrary
infinite--dimensional Lie group.

\begin{proposition}\label{pr:uniqueness}
If $\G$ is connected, there exists at most one exponentiation
of $\pi$.
\end{proposition}
\proof
Let $\rho_{i}$, $i=1,2$ be two exponentiations of $\pi$
with corresponding subspaces $\D_{\rho_{i}}$, $g\in\G$,
$p$ a smooth path in $\G$ with $p(0)=1$, $p(1)=g$ and
$\xi\in\D\subseteq\D_{\rho_{1}}\cap\D_{\rho_{2}}$. Since
\eqref{eq:local ode} determines the lift uniquely up to
multiplication by some $z\in\T$, we may assume that each
$\rho_{i}(p)$ possesses a unitary lift $\wt\rho_{i}(p)$
over $[0,1]$ satisfying \eqref{eq:local ode} and
$\wt\rho_{i}(p(0))=1$.
Set now  $F(t)=\wt\rho_{1}(p(t))\xi-\wt\rho_{2}(p(t))\xi$.
Then, $\dot{F}=\bbar{\dot{p}p^{-1}}F$ so that, by
skew--symmetry $G(t)=\|F(t)\|^{2}$ satisfies $\dot{G}
\equiv 0$ and $F(1)=F(0)=0$ whence $\wt\rho_{1}(g)=\wt
\rho_{2}(g)$ \halmos\\

\remark
If $\G$ has an exponential map and is generated by its image 
\footnote
{Such is the case for Banach--Lie groups and, for example, the
connected component of the identity of Diff$(M)$, $M$ a compact
manifold, since it is perfect \cite{Th,Ep1} and therefore
has no proper normal subgroups \cite{Ep2}.}, the above definition
of exponentiation may be weakened by requiring that
\eqref{eq:local ode} hold only for $p(t)=\exp_{\G}(tX)$. More
precisely, since continuous one--parameter groups in $PU(\H)$
lift to $U(\H)$, uniquely up to multiplication by a character of
$\IR$, $\rho(\exp_{\G}(tX))$ lifts to a one--parameter group
$\wt\rho(t)$ which may be normalised by demanding that
\begin{equation}
\left.\frac{d}{dt}\right|_{t=0}\wt\rho(t)\xi=\bbar{X}\xi
\end{equation}
It follows that $\pi(X)$ is essentially skew--adjoint
\cite[thm. VIII.10]{RS} and that
$\wt\rho(\exp_{\G}(tX))=e^{t\bbar{X}}$ so that $\rho$ is
uniquely determined by $\pi$.

\ssubsection{The exponentiation theorem}

\begin{theorem}\label{th:main}
Let $\G$ be a connected and simply connected Lie group with Lie
algebra $\L$ and $\pi:\L\longrightarrow\End(V)$ a projective
representation of $\L$ by skew--symmetric operators acting
on a dense subspace $V$ of a Hilbert space $\H$. Let
$B\in H^{2}(\L,\IR)$ be the corresponding cocyle so that
\begin{equation}
[\pi(X),\pi(Y)]=\pi([X,Y])+iB(X,Y)
\end{equation}
Assume the existence of a self--adjoint operator $A\geq 1$ on
$\H$ with
\begin{equation}
V=\bigcap_{n\geq 0}\D(A^{n})
\end{equation}
and such that, for any $n\in\IN$, $\xi\in V$ and $X\in\L$
\begin{align}
\|\pi(X)\xi\|_{n}&\leq |X|_{n+1}\|\xi\|_{n+1}
\label{eq:th assu 1}\\
\|[A,\pi(X)]\xi\|_{n}&\leq |X|_{A,n+1}\|\xi\|_{n+1}
\label{eq:th assu 2}
\end{align}
where $\|\xi\|_{n}=\|A^{n}\xi\|$ and the $|.|$ are continuous
semi--norms on $\L$. Then, $\pi$ exponentiates uniquely to $\G$.
Moreover, $\G$ leaves each $\D(A^{n})$ invariant and acts
continuously on these.
\end{theorem}

The proof of theorem \ref{th:main} depends on a number of preliminary
results. Let $I=[0,1]$, then

\begin{lemma}\label{le:integrability}
Let $X_{i}\in C^{\infty}(I^{2},\L)$, $i=1,2$ satisfy the integrability
condition
\begin{equation}\label{eq:flatness}
\pi(\partial_{1}X_{2})-\pi(\partial_{2}X_{1})=[\pi(X_{1}),\pi(X_{2})]
\end{equation}
Then, for any $\xi\in\hsmooth$,
there exists $F\in C^{\infty}(I^{2},\hsmooth)$ such that
\begin{align}
\partial_{i}F&=\pi(X_{i})F\\
       F(0,0)&=\xi
\end{align}
\end{lemma}
\proof We proceed as in the local trivialisation of flat vector
bundles. For any $(x,y)\in I$, set
\begin{equation}
F(x,y)=
\prod_{y\geq v\geq 0}\Exp(X_{2}(x,v)dv)
\prod_{x\geq u\geq 0}\Exp(X_{1}(u,0)du)\xi
\end{equation}
We claim that $F:I^{2}\rightarrow\hsmooth$ is smooth. To see this,
use the change of variable formula \eqref{eq:change of variable}
to rewrite $F(x,y)$ as $\I(Y(x,y),\I(Z(x,y),\xi))$ where
\begin{equation}
Y(x,y)(t)=yX_{2}(x,yt)
\qquad\text{and}\qquad
Z(x,y)(t)=xX_{1}(xt,0)
\end{equation}
are smooth maps $I^{2}\rightarrow C(I,\L)$ and use the chain rule
in conjunction with corollary \ref{smoothness}. By construction,
$F$ satisfies $\partial_{2}F=\pi(X_{2})F$. Set
$G=(\partial_{1}-\pi(X_{1}))F$, we
wish to show that $G=0$. This certainly holds on $I\times\{0\}$
by the very definition of $F$. Moreover, for any
$H\in C^{2}(I^{2},\hsmooth)$, \eqref{eq:flatness} implies that
$[\partial_{1}-\pi(X_{1}),\partial_{2}-\pi(X_{2})]H=0$. Thus
\begin{equation}
(\partial_{2}-\pi(X_{2}))G=
(\partial_{2}-\pi(X_{2}))(\partial_{1}-\pi(X_{1}))F=
(\partial_{1}-\pi(X_{1}))(\partial_{2}-\pi(X_{2}))F=
0
\end{equation}
so that, by uniqueness $G\equiv 0$ \halmos

\begin{lemma}\label{le:holonomy}
Let $X_{i}\in C^{\infty}(I^{2},\L)$, $i=1,2$ be such that
$[X_{1},X_{2}]=\partial_{1}X_{2}-\partial_{2}X_{1}$
and $\left.X_{2}\right|_{\{0,1\}\times I}\equiv 0$.
Then,
\begin{equation}\label{eq:holonomy}
\prod_{1\geq\tau\geq 0}\Exp(X_{1}(\tau,1)d\tau)=
e^{i\int_{0}^{1}\int_{0}^{1}B(X_{1},X_{2})dudv}
\prod_{1\geq\tau\geq 0}\Exp(X_{1}(\tau,0)d\tau)
\end{equation}
\end{lemma}
\proof
Let $\wt\L=\L\oplus c\cdot\IR$ be the extension of $\L$ by a central
element $c$ with bracket
\begin{equation}
[X\oplus tc,Y\oplus sc]=[X,Y]\oplus B(X,Y)c
\end{equation}
$\pi$ extends to a genuine representation of $\wt\L$ on $\H^{\infty}$
by letting $c$ act as multiplication by $i$, which moreover still
satisfies \eqref{eq:th assu 1}--\eqref{eq:th assu 2} provided we
set $|X\oplus tc|_{n+1}=|X|_{n+1}+|t|$ and
$|X\oplus tc|_{A,n+1}=|X|_{A,n+1}$.
Let now $Y_{i}\in C^{\infty}(I^{2},\wt\L)$ be given by $Y_{1}=X_{1}$
and 
\begin{equation}
Y_{2}(x,y)=X_{2}(x,y)+c\int_{0}^{x}B(X_{1},X_{2})(t,y)dt
\end{equation}
Then
\begin{equation}
\pi(\partial_{1}Y_{2}-\partial_{2}Y_{1})=
\pi([X_{1},X_{2}]+cB(X_{1},X_{2}))=
[\pi(X_{1}),\pi(X_{2})]=
[\pi(Y_{1}),\pi(Y_{2})]
\end{equation}
By lemma \ref{le:integrability}, we may find
$F\in C^{\infty}(I^{2},\hsmooth)$ satisfying
$\partial_{i}F=\pi(Y_{i})F$ and $F(0,0)=\xi$, an arbitrary vector
in $\hsmooth$. Since $Y_{2}(0,\cdot)\equiv 0$, we get $F(0,\cdot)=\xi$
and, by uniqueness of solutions of $\partial_{1}F=\pi(X_{1})F$,
$F(1,y)=\prod_{1\geq u\geq 0}\Exp(X_{1}(u,y)du)\xi$.
On the other hand, since $Y_{2}(1,y)=c\int_{0}^{1}B(X_{1},X_{2})(u,y)du$,
both $F(1,y)$ and
$G(y)=e^{i\int_{0}^{y}\int_{0}^{1}B(X_{1},X_{2})(u,v)dudv}F(1,0)$ are
annihilated by $(\partial_{2}-i\int_{0}^{1}B(X_{1},X_{2})(u,\cdot)du)$
whence $F(1,y)=G(y)$. Thus, \eqref{eq:holonomy} holds since both sides
coincide on $\hsmooth$ \halmos\\

\definition
For any $p\in C^{\infty}(I,\G)$, let 
$U_{p}=\prod_{1\geq\tau\geq 0}\Exp(\dot pp^{-1}(\tau)d\tau)\in U(\H)$.\\

\begin{proposition}\label{pr:Up}
The unitaries $U_{p}$ have the following properties
\begin{enumerate}
\item[(i)]
Lifting property : if $X\in\L$ defines a one--parameter group in $\G$
and $p(t)=\exp_{\G}(tX)$, then $U_{p}=e^{\bbar{X}}$.
\item[(ii)]
Invariance under reparametrisation : $U_{p\cdot\phi}=U_{p}$ for
any smooth $\phi:I\rightarrow I$ fixing the endpoints.
\item[(iii)]
Translation invariance : $U_{pg}=U_{p}$ for any $g\in\G$.
\item[(iv)]
Factorisation property : $U_{p}=U_{p_{2}}U_{p_{1}}$ where
$p_{1}(t)=p(\half{t})$ and $p_{2}(t)=p(\half{t+1})$,
\item[(v)]
Inversion property :
$U_{p}^{*}=U_{\check p}$ where $\check{p}(t)=p(1-t)$.
\item[(vi)]
Projective homotopy invariance :
If $p_{i}$, $i=0,1$ are homotopic relative to their endpoints
through $H\in C^{\infty}(I^{2},\G)$, then
\begin{equation}\label{change of path}
U_{p_{1}}=e^{i\int_{H(I^{2})}B}U_{p_{0}}
\end{equation}
where
$\int_{H(I^{2})}B=
 \int_{0}^{1}\int_{0}^{1}
 B(\partial_{1}H\cdot H^{-1},\partial_{2}H\cdot H^{-1})dx_{1}dx_{2}$.
\end{enumerate}
\end{proposition}
\proof
(i) Since $\dot pp^{-1}\equiv X$, the claim follows from the definition
of product integrals.\\
(ii) We have
$\frac{d}{dt}p\cdot\phi(p\cdot\phi)^{-1}=\phi'\dot pp^{-1}\cdot\phi$
and the claim follows from the change of variable formula
\eqref{eq:change of variable}.\\
(iii) follows from $\dot{(pg)}(pg)^{-1}=\dot pp^{-1}$.\\
(iv) follows from the semigroup property of product integrals.\\
(v) follows from the inversion formula \eqref{eq:inversion}.\\
(vi)
Let $X_{i}=\partial_{i}H\cdot H^{-1}\in C^{\infty}(I^{2},\L)$ so that
$X_{2}(0,\cdot)=X_{2}(1,\cdot)\equiv 0$. Since $\partial_{i}H=X_{i}H$,
we have
\begin{equation}
\partial_{1}\partial_{2}H=
\partial_{1}X_{2}H+X_{2}\partial_{1}H=
\partial_{1}X_{2}H+X_{2}X_{1}H
\end{equation}
and similarly
$\partial_{2}\partial_{1}H=\partial_{2}X_{1}H+X_{1}X_{2}H$.
Substracting and multiplying by $H^{-1}$ to the right we get
\begin{equation}
\partial_{1}X_{2}-\partial_{2}X_{1}=[X_{1},X_{2}]
\end{equation}
The result now follows from lemma \ref{le:holonomy} \halmos\\

{\sc Proof of theorem \ref{th:main}}.
Define a map $\rho:\G\longrightarrow PU(\H)$ (resp. $U(\H)$ if
$B=0$) as follows. Let $g\in\G$, $p:I\rightarrow\G$ a smooth
path with $p(0)=1$ and $p(1)=g$ and set $\rho(g)=U_{p}$. By
\eqref{change of path}, this is a well--defined map. To check
that $\rho$ is a homomorphism, let $g,h\in\G$ and
$p\in C^{\infty}(I,\G)$ such that $p(0)=1$, $p(\half{1})=g$,
$p(1)=hg$. By the semigroup property
\begin{equation}
U_{p}=
\prod_{1\geq\tau\geq 0}\Exp(\dot pp^{-1}(\tau)d\tau)=
\prod_{1\geq\tau\geq \half{1}}\Exp(\dot pp^{-1}(\tau)d\tau)
\prod_{\half{1}\geq\tau\geq 0}\Exp(\dot pp^{-1}(\tau)d\tau)=
U_{p_{2}}U_{p_{1}}
\end{equation}
where $p_{1}(t)=p(2t)$ goes from $1$ to $g$ and $p_{2}(t)=p(2t-1)$
goes from $g$ to $hg$. By translation invariance,
$U_{p_{2}}=U_{p_{2}g^{-1}}=\rho(h)$
and therefore $\rho(hg)=U_{p}=\rho(h)\rho(g)$. The strong
continuity of $\rho$ follows from corollary \ref{co:soI cont}.
Set $\D_{\rho}=\H^{1}$, a $\G$--invariant subspace by
\S \ref{ss:volterra}. To check that the differential of $\rho$
on $\D_{\rho}$ is $\pi$, fix $p\in C^{\infty}(\IR,\G)$ with
$p(0)=1$. Then, \begin{equation}
\wt\rho(p(t))=\prod_{t\geq\tau\geq 0}\Exp(\dot{p}p^{-1}(\tau)d\tau)
\end{equation}
is a lift of $\rho(p(t))$ (resp. equals $\rho(p(t))$ if $B=0$)
which, by theorem \ref{th:homo eq} satisfies
\begin{equation}
\frac{d}{dt}\wt\rho(p)\xi=\bbar{\dot{p}p^{-1}}\wt\rho(p)\xi
\end{equation}
for any $\xi\in\D_{\rho}$. Finally, $\G$ leaves each $\D(A^{n})$
invariant by \S \ref{ss:volterra} and acts continuously on them
by corollary \ref{co:soI cont} \halmos\\

\remark
The main conclusion of theorem \ref{th:main} depends only on the fact
that $\L$ is a Lie algebra and not on the fact that it has an underlying
Lie group. More precisely, to any topological Lie algebra $\L$, one
can associate an abstract group $T(\L)$ called the Thurston group of
$\L$ in the following way \cite[5.5]{Mi2}. The elements of $T(\L)$
are equivalence 
classes of smooth paths $I\rightarrow\L$ for the relation
$u_{0}\sim u_{1}$ if there exist $X_{i}\in C^{\infty}(I^{2},\L)$,
$i=1,2$ satisfying $X_{1}(t,0)=u_{0}(t)$, $X_{1}(t,1)=u_{1}(t)$,
$\left.X_{2}\right|_{\{0,1\}\times I}\equiv 0$ as well as the
integrability condition
\begin{equation}
[X_{1},X_{2}]=\partial_{1}X_{2}-\partial_{2}X_{1}
\end{equation}
Any connected and simply connected Lie group $\G$ maps homomorphically
into the Thurston group of its Lie algebra by associating to $g\in\G$
the class of $\dot pp^{-1}$ where $p$ is any smooth path in $\G$ with
$p(0)=1$, $p(1)=g$ and this map is an isomorphism if $\G$ is regular
in the sense explained in the introduction.
Now if $\pi$ is a projective representation of $\L$ satisfying the
assumptions of theorem \ref{th:main}, lemma \ref{le:holonomy} shows
that $\pi$ yields a (projective) unitary representation of $T(\L)$
and therefore one of the Lie group underlying $\L$ if one such exists.

\ssubsection{Smoothness of central extensions of
$\G$ arising from exponentiated representations}
\label{ss:smooth ext}

A projective unitary representation $\rho:\G\longrightarrow PU(\H)$
lifts to a unitary representation of the continuous central extension
$\pback$ of $\G$ obtained by pulling back the canonical central
extension
\begin{equation}\label{eq:can ext}
1\rightarrow\T\rightarrow U(\H)\xrightarrow{p} PU(\H)\rightarrow 1
\end{equation}
to $\G$. Explicitly, 
\begin{equation}\label{eq:pback}
\pback=\{(g,V)\in\G\times U(\H)|\rho(g)=p(V)\}
\end{equation}
acts on $\H$ by $(g,V)\xi=V\xi$. For classification purposes, it is
often useful to know that $\pback$ is a Lie group. This is so if $\G$
is finite--dimensional since any local continuous cocycle of $\pback$
may be regularised within its cohomology class by convolving it with
a smooth function on $\G\times\G$ (see \eg \cite[lemma 7.20]{Va}). In
the absence of a Haar measure this ceases to be obvious but continues
to hold for the class of representations considered in this paper.

\begin{proposition}\label{pr:smooth ext}
Let $\pi:\L\longrightarrow\End(V)$ be a projective representation
of $\L$ satisfying the assumptions of theorem \ref{th:main}, $B$
its cocycle and $\rho:\G\longrightarrow PU(\H)$ its exponentiation.
Then, $\pback$ is a smooth central extension of $\G$. In particular,
its isomorphism class is uniquely determined by the Lie algebra
cocycle of $\pback$ which is equal to $B$.
\end{proposition}
\proof
It is sufficient to exhibit a local trivialisation of $\pback$ the
corresponding local multiplication and inversion of which are smooth.
We begin by trivialising \eqref{eq:can ext} as in \cite{Ba}. Fix
$\xi\in\H$ of norm 1 and consider the open set
\begin{equation}
U_{\xi}=\{[u]\in PU(\H)|\thinspace |(u\xi,\xi)|>0\}
\end{equation}
where $[u]$ is the equivalence class of $u\in U(\H)$ in $PU(\H)$.
Define a function
\begin{equation}
\phase:p^{-1}(U_{\xi})\longrightarrow\T,
\quad
u\longrightarrow\phase(u)=\frac{(u\xi,\xi)}{|(u\xi,\xi)|}
\end{equation}
and notice that $\phase(e^{i\theta}u)=e^{i\theta}\phase(u)$
so that the map $\phi:p^{-1}(U_{\xi})\longrightarrow U_{\xi}\times\T$,
$\phi(u)=([u],\phase(u))$
is a $\T$--equivariant local trivialisation with inverse
$\phi^{-1}([u],z)=u\phase(u)^{-1}z$. The corresponding local
multiplication and inversion on $\G\times\T$, namely
\begin{align}
x\star y&=\phi(\phi^{-1}x\cdot\phi^{-1}y)\\
i(x) &=\phi\bigl((\phi^{-1}x)^{*}\bigr)
\end{align}
read, using $\phase(u^{*})=\overline{\phase(u)}$,
\begin{align}
([u],z)\star([v],w)&=\Bigl([uv],zw
\frac{\phase(uv)}{\phase(u)\phase(v)}\Bigr)\\
i([u],z)&=([u^{*}],z^{-1})
\end{align}
where the quotient
$\phase(uv)\phase(u)^{-1}\phase(v)^{-1}$
is independent of the choice of the lifts $u,v$ of $[u],[v]\in PU(\H)$.
Pulling back by $\rho$, we obtain a local trivialisation of $\pback$
with group laws
\begin{align}
(g,z)\star(h,w)&=
\Bigl(gh,zw\frac
{\phase(\rho(g)\rho(h))}{\phase(\rho(g))\phase(\rho(h))}
\Bigr)\label{eq:local mult}\\
i(g,z)&=(g^{-1},z^{-1})\label{eq:local inv}
\end{align}
We claim that if $\xi\in\H^{\infty}$, the local multiplication
\eqref{eq:local mult} is smooth. To see this, let $e:\L\longrightarrow\G$
be a local chart mapping $0$ to $1$. $e$ defines a local smooth
embedding $p:\G\longrightarrow C^{\infty}(I,\G)$ mapping $g\in\G$
to the path $p(g,t)=e(te^{-1}(g))$ with endpoints $1$, $g$ and,
by the construction of $\rho$ given in theorem \ref{th:main},
\begin{equation}\label{eq:a lift}
\wt\rho:
g\longrightarrow
\prod_{1\geq\tau\geq 0}\Exp(\dot p(g,\tau)\cdot p(g,\tau)^{-1}d\tau)
\end{equation}
is a local lift of $\rho$. By corollary \ref{smoothness}, the map
$\G\times\H^{\infty}\longrightarrow\H^{\infty}$,
$(g,\xi)\longrightarrow\wt\rho(g)\xi$ is smooth and therefore
so is the local multiplication \eqref{eq:local mult}.\\

We now compute the local adjoint action and Lie bracket of $\pback$.
It will be more convenient to assume that the derivative of $e$ at
$0$ is the identity, which may be achieved by pre--multiplying $e$
by $(D_{0}e)^{-1}$. Let $X\in\L$ and $x\in\IR$, then, identifying
the Lie algebra of $\T$ with $i\IR$, we find
\begin{equation}\label{eq:local Ad}
\begin{split}
\Ad(g,z)\medspace X\oplus ix
&=
\diff{t}{0}(g,z)\star(e(tX),e^{itx})\star(g,z)^{-1}\\
&=
\diff{t}{0}\left(ge(tX)g^{-1},
\frac
{\phase(\wt\rho(g)\wt\rho(e(tX))\wt\rho(g)^{*})}
{\phase(\wt\rho(e(tX)))}
e^{itx}\right)
\end{split}
\end{equation}
Set $p(t)=e(tX)$, then
\begin{equation}
\wt\rho(e(tX))=
\prod_{1\geq\tau\geq 0}\Exp(\tau\dot pp^{-1}(\tau t)d\tau)=
\prod_{t\geq\tau\geq 0}\Exp(\dot pp^{-1}(\tau)d\tau)
\end{equation}
by the change of variable formula \eqref{eq:change of variable}
and it follows by theorem \ref{th:homo eq} that for any
$\eta\in\H^{\infty}$,
\begin{equation}\label{eq:stanco}
\diff{t}{0}\wt\rho(e(tX))\eta=
\pi(X)\eta
\end{equation}
Thus
\begin{equation}
\diff{t}{0}\phase(\wt\rho(e(tX)))=
\diff{t}{0}\frac
{(\wt\rho(e(tX))\xi,\xi)}
{\sqrt{(\wt\rho(e(tX))\xi,\xi)(\xi,\wt\rho(e(tX))\xi)}}=
(\pi(X)\xi,\xi)
\end{equation}
and similarly,
\begin{equation}
 \diff{t}{0}\phase(\wt\rho(g)\wt\rho(e(tX))\wt\rho(g)^{*})=
(\wt\rho(g)\pi(X)\wt\rho(g)^{*}\xi,\xi)
\end{equation}
since $\G$ leaves $\H^{\infty}$ invariant. The local adjoint
action \eqref{eq:local Ad} is therefore given by
\begin{equation}
\Ad(g,z)\medspace X\oplus ix=
gXg^{-1}\oplus
\Bigl(ix+(\wt\rho(g)\pi(X)\wt\rho(g)^{*}\xi,\xi)-(\pi(X)\xi,\xi)\Bigr)
\end{equation}

Take now $g=e(sY)$, $Y\in\L$. By \eqref{eq:stanco} and the
skew--symmetry of $\pi(Y)$,
\begin{equation}
\diff{s}{0}\wt\rho(e(sY))^{*}\xi=-\pi(Y)\xi
\end{equation}
so that the Lie bracket on $\L\oplus i\IR$ is
\begin{equation}
\begin{split}
[Y\oplus iy,X\oplus ix]
&=\diff{s}{0}\Ad(e(sY),e^{isy})\medspace X\oplus ix\\[1ex]
&=[Y,X]\oplus ([\pi(Y),\pi(X)]\xi,\xi)\\[.8ex]
&=[Y,X]\oplus \Bigl(iB(Y,X)+(\pi([Y,X])\xi,\xi)\Bigr)
\end{split}
\end{equation}
Thus, the Lie algebra cocycle of $\pback$ is $B(Y,X)-i(\pi([Y,X])\xi,\xi)$
which is cohomologous to $B$. Finally, the fact that a smooth
central extension of $\G$ is uniquely determined by its Lie algebra cocycle
is proved in \cite[page 54]{PS} \footnote{the proof is only given for loop
groups but works {\it verbatim} for any connected and simply connected Lie
group.} \halmos\\

\remark Proposition \ref{pr:smooth ext} may also be proved in the
following way. By \eqref{change of path}, $B$, when regarded as a
right--invariant, closed two--form on $\G$ is integral, \ie its value
on closed two--cycles is an integral multiple of $2\pi$. It follows
that there exists a smooth central extension $\wt\G$ of $\G$ by $\T$
with Lie algebra cocycle $B$ which may be described as follows
\cite[prop. 4.4.2]{PS}. Let $\P\G$ be the space of piecewise smooth
paths $I\rightarrow\G$ with $p(0)=1$. The concatenation of pointed
paths defined by
\begin{equation}
p\vee q(t)=
\left\{\begin{array}{rcl}
q(2t)&\text{if}&0\leq t\leq\half{1}\\
p(2t-1)q(1)&\text{if}&\half{1}\leq t\leq 1
\end{array}\right.
\end{equation}
induces a monoidal structure on $\P\G\times\T$ and $\wt\G$ is the
quotient of $\P\G\times\T$ by the equivalence relation
\begin{equation}
(p,z)\sim(q,w)
\quad\Longleftrightarrow\quad
p(1)=q(1)
\thickspace\thickspace\text{and}\thickspace\thickspace
e^{i\int_{\sigma}\beta}=w\overline{z}
\end{equation}
where $\sigma$ is any two--cycle with boundary $p\vee\check{q}$
and $\check{q}(t)=q(1-t)q(1)^{-1}$. The construction of the
exponentiation $\rho$ of $\pi$ given in theorem \ref{th:main}
then shows that the map
\begin{equation}
\P\G\times\T\longrightarrow\G\times U(\H),
\quad
(p,z)\longrightarrow(p(1),zU_{p})
\end{equation}
descends to an isomorphism of central extensions $\wt\G\cong\pback$
and in particular that $\pback$ is smooth.\\

\remark By proposition \ref{pr:smooth ext} and \cite[prop. 7.1]{Se},
the topological type of $\pback$ as a principal circle bundle over
$\G$ is determined by the image of $B$ in $H^{2}(\G,\IR)$.
In particular, if $B=dA$ for some 1--form $A$ on $\G$, then,
by \eqref{change of path}
\begin{equation}
U_{p}e^{i\int_{p}A}=U_{q}e^{i\int_{q}A}
\end{equation}
for any two homotopic paths $p,q$ in $\G$, so that the map
$g\longrightarrow (g,U_{p}e^{i\int_{p}A})$ where
$p(0)=1$, $p(1)=g$ is a section of $\pback$.

\ssection{Applications}
\label{se:applications}

\ssubsection
{Positive energy representations of Diff$\mathbf{(S^{1})}$ and loop groups}

Let $\Diffp$ be the group of orientation--preserving diffeomorphisms
of $S^{1}$. It is isomorphic to the quotient of the group $\D$
of diffeomorphisms $\phi$ of $\IR$ satisfying
$\phi(x+2\pi)=\phi(x)+2\pi$ by the subgroup of translations by
multiples of $2\pi$. Since $\D$ is contractible through the map
$(\phi,t)\rightarrow t\id+(1-t)\phi$, $\Diffp$ is connected and
$\D$ is the universal covering group of $\Diffp$. The Lie algebra
of $\Diffp$ and $\D$ is $\vect$, the smooth real vector fields on
$S^{1}$, with bracket
\begin{equation}\label{eq:vect bracket}
[f\der{\theta},g\der{\theta}]=
(f'g-fg')\der\theta
\end{equation}
The Virasoro algebra $\Vir$ is by definition the central extension
of the Lie algebra $\vectpolc$ of complex vector fields with finite
Fourier series with respect to the cocycle
\begin{equation}\label{eq:vect cocycle}
\omega(f\der{\theta},g\der{\theta})=
\frac{1}{12}\int_{0}^{2\pi}(f''+f)g'\frac{d\theta}{2\pi}
\end{equation}
It is spanned by $L_{n}=-ie^{in\theta}\frac{d}{d\theta}$ , $n\in\IZ$
and a central element $\kappa$ in terms of which the bracket reads
\begin{equation}\label{eq:vir comm}
[L_{m},L_{n}]=(m-n)L_{m+n}+\delta_{m+n,0}\frac{m^{3}-m}{12}\kappa
\end{equation}
A highest weight representation of $\Vir$ is a representation $V$
such that $\kappa$ acts as multiplication by a scalar $c$ and $V$
is generated over the enveloping algebra of the $L_{n}$, $n<0$ by
an $L_{0}$--eigenvector $\Omega$ annihilated by the $L_{n}$, $n>0$.
If $L_{0}\Omega=h\Omega$, it follows from \eqref{eq:vir comm} that
$L_{0}$ is diagonal with finite--dimensional eigenspaces and that
its spectrum is contained in $h+\IN$. The pair $(c,h)$ is called
the highest weight of $V$. Let $\overline{\medspace\cdot\medspace}$
be the anti--linear anti--automorphism of $\Vir$ acting as $-1$ on
real vector fields so that $\overline{L_{n}}=L_{-n}$ and, by
\eqref{eq:vir comm}, $\overline{\kappa}=\kappa$. $V$ is called
unitarisable if it possesses an inner product $(\cdot,\cdot)$
such that $(X\xi,\eta)=(\xi,\overline{X}\eta)$ for any $X\in\Vir$
and $\xi,\eta\in V$. In that case,
$V$ is irreducible and $c,h\in\IR$. In fact $c,h\geq 0$ since,
for any $n\geq 1$, \eqref{eq:vir comm} yields
\begin{equation}
0\leq
(L_{-n}\Omega,L_{-n}\Omega)=
2nh+\frac{n^{3}-n}{12}c
\end{equation}
so that $L_{0}$ has non--negative spectrum. The values
of $(c,h)$ for which there exists a unitarisable module with highest
weight $(c,h)$ have been classified by the joint results of
Friedan--Qiu--Shenker \cite{FQS1,FQS2} and Goddard--Kent--Olive
\cite{GKO}. They are $\{(c,h)|c\geq 1,\medspace h\geq 0\}$
together with the discrete series contained in the region
$0\leq c<1$, $h\geq 0$ and parametrised by
\begin{align}
c(m)      &=1-\frac{6}{(m+2)(m+3)}\\
h_{p,q}(m)&=\frac{((m+3)p-(m+2)q)^{2}-1}{4(m+2)(m+3)}
\end{align}
where $m\geq 1$, $p=1\ldots m+1$ and $q=1\ldots p$. We won't need
however to rely on this classification. The following result was
conjectured by Kac and proved in special cases by Segal \cite{Se}
and Neretin \cite{Ner} and in the general case by Goodman and Wallach
\cite[thm. 4.2]{GoWa2}

\begin{theorem}
Let $(\pi,V)$ be a unitarisable highest weight representation of the
Virasoro algebra. Then, $\pi$ exponentiates uniquely to a projective
unitary representation of $\Diffp$ on the Hilbert space completion
$\H$ of $V$.
\end{theorem}
\proof
Let $(c,h)$ be the highest weight of $V$ and $A=1+\bbar{L_{0}}$
acting on $\H$. Using simple $\mathfrak{sl}_{2}(\IC)$ arguments
relying on the positivity of the spectrum of $\pi(L_{0})$, Goodman
and Wallach showed  \cite[prop. 2.1]{GoWa2} that for any $t\in\IR$,
$\xi\in V$ and $X=\sum_{n}a_{n}e^{in\theta}\der{\theta}\in\vectpolc$
\begin{equation}
\|\pi(X)\xi\|_{t}\leq 
2^{\half{1}}\|X\|_{|t|}\|\xi\|_{t+1}+
M\|X\|_{|t|+1}\|\xi\|_{t+\half{1}}+
M\|X\|_{|t|+\half{3}}\|\xi\|_{t}
\end{equation}
where $M=(c/12)^{\half{1}}$, $\|\xi\|_{t}=\|A^{t}\xi\|$ and
$\|X\|_{s}=\sum_{n}(1+|n|)^{s}|a_{n}|$.
Thus, $\pi$ extends to a projective unitary representation of
$\vect$ on the space of smooth vectors of $A$ satisfying
\eqref{eq:th assu 1}--\eqref{eq:th assu 2} where
\begin{align}
|X|_{n+1}&=2^{\half{1}}\|X\|_{n}+M(\|X\|_{n+1}+\|X\|_{n+\half{3}})\\
|X|_{A,n+1}&=|[L_{0},X]|_{n+1}
\end{align}
By theorem \ref{th:main}, $\pi$ exponentiates to a projective unitary
representation $\rho$ of the universal cover $\D$ of $\Diffp$. We claim
that the kernel of the covering map, \ie translations
$T_{2\pi n}=\exp_{\D}(2\pi inL_{0})$ by multiples of $2\pi$, acts by
scalars. Indeed, since the spectrum of $\pi(L_{0})$ is contained in
$h+\IN$, we have, by the definition of $\rho$ and (i) of proposition
\ref{pr:Up},
\begin{equation}
\rho(T_{2\pi n})=e^{2\pi n i\bbar{L_{0}}}=e^{2\pi inh}
\end{equation}
and the projective action of $\D$ factors through $\Diffp$
\halmos\\

Let now $G$ be a compact, connected and simply connected simple
Lie group and $LG=C^{\infty}(S^{1},G)$ its loop group. Since
$LG=\Omega G\rtimes G$ where $\Omega G$ is the space of based
loops and $G$ that of constant ones, we get
$\pi_{0}(LG)=\pi_{1}(G)=0$ and
$\pi_{1}(LG)=\pi_{2}(LG)\oplus\pi_{1}(G)=0$
so that $LG$ is connected and simply connected. The Lie algebra
of $LG$ is $L\g=C^{\infty}(S^{1},\g)$ where $\g$ is the Lie algebra
of $G$ and has complexification $L\gc$. $L\g$ has a distinguished
cocycle which generates $H^{2}(L\g,\IR)$ \cite[prop. 4.2.4]{PS},
namely
\begin{equation}
B(X,Y)=\int_{0}^{2\pi}\<X,Y'\>\frac{d\theta}{2\pi}
\end{equation}
where $\<\cdot,\cdot\>$ is the basic inner product, \ie the
multiple of the Killing form for which the highest root $\theta$
has squared length 2. $L\g$ has a dense subalgebra $\lpol\g$
consisting of all $\g$--valued trigonometric polynomials
and the central extension of $\lpol\gc$ corresponding to $B$
is usually denoted by $\wt\gc$. It is spanned by elements 
$x(n)=x\otimes e^{in\theta}$, $x\in\gc$, $n\in\IZ$ and a
central element $\kk$ with bracket
\begin{equation}
[x(m),y(n)]=[x,y](m+n)+m\delta_{m+n,0}\<x,y\>\kk
\end{equation}

Since $B$ is invariant under the action of $\Diffp$ on $L\g$
by reparametrisation, one may form the semi--direct product
$\wt\gc\rtimes\vectpolc$. By definition, the affine Kac--Moody
algebra $\wh{\gc}$ is the subalgebra $\wt\gc\rtimes\IC\cdot L_{0}$.
Let $T$ be a maximal torus in $G$ with Lie algebra $\t$.
A highest weight representation of $\wh\gc$ is a representation
where the central element $\kk$ acts by a scalar and which is
generated over the enveloping algebra of the $x(n)$, $n<0$ or
$n=0$ and $x$ a root vector corresponding to a negative
root, by a vector $\Omega$ diagonalising the action of
$\tc\times\IC\cdot L_{0}$ and annihilated by the $x(n)$, $n>0$
or $n=0$ and $x$ a root vector corresponding to a positive root.
Thus, for any $t\in\tc$
\begin{align}
  \kk\Omega &= \ell\Omega\\
    t\Omega &= \lambda(t)\Omega\\
L_{0}\Omega &= h\Omega
\end{align}
for some $\ell,h\in\IC$ and $\lambda\in\tc^{*}$. In particular,
since $[L_{0},x_{n}]=-nx(n)$, $L_{0}$ is diagonal with
finite--dimensional eigenspaces and spectrum contained in
$h+\IN$.\\

A fundamental feature of highest weight representations is that
the action of $\wh\gc$ extends to one of $\wt\gc\rtimes\Vir$,
provided $\ell+h^{\vee}\neq 0$ where $h^{\vee}$ is the dual Coxeter
number of $\gc$. This is obtained via the Segal--Sugawara formulae
by letting $L_{n}$, $n\neq 0$ act as
\begin{equation}\label{eq:ss n}
\frac{1}{2(\ell+h^{\vee})}
\sum_{\substack{i\\m\in\IZ}}x_{i}(-m)x^{i}(m+n)
\end{equation}
$L_{0}$ as
\footnote
{It is easy to see that \eqref{eq:ss 0} differs from the original
action of $L_{0}$ by an additive constant equal to the difference
of their lowest eigenvalues, namely
$h-C_{\lambda}/2(\ell+h^{\vee})$ where $C_{\lambda}$ is the Casimir
of the irreducible $G$--module with highest weight $\lambda$.}
\begin{equation}\label{eq:ss 0}
\frac{1}{2(\ell+h^{\vee})}\left(
\sum_{i}x_{i}(0)x^{i}(0)+
2\sum_{\substack{i\\n>0}}x_{i}(-n)x^{i}(n)\right)
\end{equation}
and $\kappa$ as multiplication by $\dim(G)\ell/(\ell+h^{\vee})$
where $x_{i},x^{i}$ are dual basis of $\gc$ for the basic inner
product \cite[\S 9.4]{PS}, \cite[thm. 10.1]{KR}.
Let $\overline{\medspace\cdot\medspace}$ be the anti--linear
anti--automorphism on $\wh\gc$ acting as -1 on $\lpol\g$,
$\overline{\kk}=\kk$ and $\overline{L_{0}}=L_{0}$ and define
unitarisable highest weight representations $V$ of $\wh\gc$
accordingly. If $(\ell,\lambda)$ is the highest weight of
$V$ this is equivalent to requiring that $V$ be irreducible
and that $\ell\in\IN$, $h\in\IR$, $\lambda$ is an integral
dominant weight of $\G$ and $\<\lambda,\theta\>\leq\ell$
\cite[thm. 11.7]{K}. In that case, the action of $\Vir$ 
given by \eqref{eq:ss n}--\eqref{eq:ss 0} is also unitarisable
and $V$ splits into an orthogonal direct sum of highest weight
representations $V_{i}$ of $\Vir$ of highest weights
$(\ell\dim(G)/(\ell+h^{\vee}),h_{i})$.
The following theorem was first proved in special cases by
Segal \cite{Se} and in the general case by Goodman and Wallach
\cite[thm 6.7]{GoWa1}

\begin{theorem}
Let $(\pi,V)$ be an integrable highest weight representation
of $\wh\gc$. Then, $\pi$ exponentiates uniquely to a projective
unitary representation of $LG$ on the Hilbert space completion
$\H$ of $V$ extending to $LG\rtimes\Diffp$.
\end{theorem}
\proof
Let $(\ell,\lambda)$ be the highest weight of $V$ and extend
$\pi$ to a representation of $\wt\gc\rtimes\Vir$ by
\eqref{eq:ss n}--\eqref{eq:ss 0}. Let $A=1+\bbar{L_{0}}$ acting
on $\H$. As noted by Goodman and Wallach, the Segal--Sugawara
formula \eqref{eq:ss 0} for $L_{0}$ readily implies that for
any $\xi\in V$, $X=\sum_{n}a_{n}e^{in\theta}\in\lpol\gc$,
$f\der{\theta}=\sum_{n}b_{n}e^{in\theta}\der{\theta}\in\vectpolc$
and $t\in\IR$
\begin{align}
\|\pi(X)\xi\|_{t}&\leq
(\ell+1)\|X\|_{|t|+\half{1}}\|v\|_{t+\half{1}}\\
\|\pi(f\der{\theta})\xi\|_{t}&\leq
\dim(G)\|f\der{\theta}\|_{|t|+\half{3}}\|\xi\|_{t+1}
\end{align}
where
$\|\xi\|_{s}=\|A^{s}\xi\|$,
$\|X\|_{s}=\sum_{n}(1+|n|)^{s}\|a_{n}\|$ and
$\|f\der{\theta}\|_{s}=\sum_{n}(1+|n|)^{s}\|b_{n}\|$
\cite[lemmas 3.2,3.3]{GoWa1}. By continuity, $\pi$ extends
to a projective unitary representation of $L\g\rtimes\vect$
on the space of smooth vectors of $A$ satisfying
\eqref{eq:th assu 1}--\eqref{eq:th assu 2} where
\begin{align}
|X\oplus f\der{\theta}|_{n+1}&=
(\ell+1)\|X\|_{n+\half{1}}+
\dim(G)\|f\der{\theta}\|_{n+\half{3}}\\
|X\oplus f\der{\theta}|_{A,n+1}&=
|[L_{0},X\oplus f\der{\theta}]|_{n+1}=
|X'\oplus f'\der{\theta}|_{n+1}
\end{align}
By theorem \ref{th:main}, $\pi$ exponentiates uniquely to a projective
unitary representation $\rho$ of $LG\rtimes\D$. Since the spectrum of
$L_{0}$ is contained in $C_{\lambda}/2(\ell+h^{\vee})+\IN$,
\begin{equation}
\rho(\exp_{\D}(2\pi inL_{0}))=
e^{2\pi in\bbar{L_{0}}}=
e^{2\pi iC_{\lambda}/2(\ell+h^{\vee})}
\end{equation}
and $\rho$ factor through a representation of $LG\rtimes\Diffp$ \halmos

\ssubsection{Unitary representations of finite--dimensional Lie algebras}

We now derive Nelson's exponentiation criterion from theorem
\ref{th:main}. We shall need the following result which is
stated as part of lemma 5.2 in \cite{Ne1} but not fully
proved there. I am grateful to Professor Z. Magyar for
showing me how to complete its proof.

\begin{lemma}\label{le:magyar}
Let $\pi:\g\longrightarrow\End(V)$ be a representation of a
finite--dimensional Lie algebra by skew--symmetric operators
on a dense subspace of a Hilbert space $\H$, $X_{i}$ a basis
of $\g$ and $\Delta=-\sum_{i}\pi(X_{i})^{2}$ the corresponding
Laplacian. If $\Delta$ is essentially self--adjoint, the
closures of the operators $\pi(X_{i})$ leave
$\H^{\infty}=\bigcap_{n\geq 0}D(\overline{\Delta}^{n})$
invariant and
\begin{equation}\label{eq:laplacian}
\left.\overline{\Delta}\right|_{\H^{\infty}}=
\sum_{i}\bbar{X_{i}}^{2}
\end{equation}
\end{lemma}
\proof
By skew--symmetry of the $\pi(X_{i})$ and self--adjointness
of $\overline{\Delta}$
\begin{equation}\label{eq:inclusion}
\sum_{i}\bbar{X_{i}}^{2}\subseteq
\sum_{i}{\pi(X_{i})^{*}}^{2}\subseteq
\Delta^{*}=\overline{\Delta}
\end{equation}
Thus, if
\begin{equation}
\K^{\infty}=\bigcap_{\substack{n\geq 0\\i_{1}\ldots i_{n}}}
\D(\bbar{X_{i_{1}}}\cdots\bbar{X_{i_{n}}})
\end{equation}
then $\K^{\infty}$ is invariant under the $\bbar{X_{i}}$
and, by \eqref{eq:inclusion} under $\overline{\Delta}$ so that
$\K^{\infty}\subseteq\H^{\infty}$. The converse 
inclusion is proved in \cite{Ne1} as formula (5.4) of lemma 5.2.
It follows that $\H^{\infty}=\K^{\infty}$ is invariant under
the $\bbar{X_{i}}$ \halmos

\begin{theorem}[Nelson]
Let $\pi:\g\longrightarrow\End(\D)$ be a representation of a
finite--dimensional Lie algebra by skew--symmetric operators
on a dense subspace of a Hilbert space $\H$, $X_{i}$ a basis
of $\g$ and $\Delta=\sum_{i}\pi(X_{i})^{2}$ the corresponding
Laplacian. If $\Delta$ is essentially self--adjoint, then $\pi$
exponentiates uniquely to a unitary representation of the
underlying connected and simply connected lie group $G$.
\end{theorem}
\proof
Let $\H^{\infty}=\bigcap_{n\geq 0}\D(\overline{\Delta}^{n})$ be
as in lemma \ref{le:magyar}. For any $X=\sum c_{i}X_{i}\in\g$, set
$\wt\pi(X)=
 \sum c_{i}\left.\bbar{X_{i}}\right|_{\H^{\infty}}
 \in\End(\H^{\infty})$
so that $\wt\pi$ is a linear action by skew--symmetric operators
and extends $\pi$. $\wt\pi$ is a representation since if
$Z=[X,Y]\in\g$, then for any $\xi\in\H^{\infty}$ and $\eta\in\D$,
\begin{equation}
\begin{split}
([\wt\pi(X),\wt\pi(Y)]\xi,\eta)
&= -(\xi,[\pi(X),\pi(Y)]\eta)\\
&= -(\xi,\pi(Z)\eta)\\
&=  (\wt\pi(Z)\xi,\eta)
\end{split}
\end{equation}
Set now $A=1-\overline{\Delta}$. By \eqref{eq:laplacian}, the
Laplacian of $\wt\pi$ is the restriction of $\overline{\Delta}$
to $\H^{\infty}$ so that the estimates
\eqref{eq:th assu 1}--\eqref{eq:th assu 2} follow
from lemma 6.3 of \cite{Ne1} and $\wt\pi$ exponentiates to $G$
by theorem \ref{th:main} \halmos

\end{document}